\documentclass[12pt,twoside]{amsart}
\usepackage{amssymb}
\usepackage{amsmath}

\newtheorem{theorem}{Theorem}[section]
\newtheorem{lemma}[theorem]{Lemma}
\newtheorem{proposition}[theorem]{Proposition}
\newtheorem{corollary}[theorem]{Corollary}
\newtheorem{remark}[theorem]{Remark}
\theoremstyle{definition}
\newtheorem{definition}[theorem]{Definition}
\newtheorem{example}[theorem]{Example}

\newcommand {\zz}{\mathbb{Z}}
\newcommand {\nn}{\mathbb{N}}
\newcommand {\qq}{\mathbb{Q}}

\newcommand{\cl}{{\mathcal L}}
\newcommand{\oo}{{\mathcal O}}
\newcommand{\ii}{{\mathcal I}}
\newcommand{\ns}{{\mathcal N}}
\newcommand {\pp}{\mathbb{P}}
\renewcommand {\aa}{\mathbb{A}}

\newcommand{\Hom}{\operatorname{Hom}}
\newcommand{\homs}{\mathcal{H}om}

\newcommand{\Proj}{\operatorname{Proj}}

\newcommand{\coker}{\operatorname{coker}}
\newcommand{\lra}{\longrightarrow}

\newcommand{\Hilb}{\mathcal{H}ilb}

\title[Doubling rational normal curves]{Doubling rational normal curves}

\author[Roberto Notari, Ignacio Ojeda, Maria
Luisa Spreafico]{Roberto Notari, Ignacio Ojeda, Maria Luisa
Spreafico}

\address{Dipartimento di Matematica, Politecnico di Milano, I-20133 Milano,
Italy}
\email{roberto.notari@polimi.it}

\address{Departamento de Matem\'aticas, Universidad de Extremadura,
E-06071 Badajoz, Espa\~na.} \email{ojedamc@unex.es}

\address{Dipartimento di Matematica, Politecnico di Torino, I-10129 Torino, Italy}
\email{maria.spreafico@polito.it}

\date{}

\begin{document}

\begin{abstract} In this paper, we study double structures supported on rational
normal curves. After recalling the general construction of double
structures supported on a smooth curve described in \cite{fer}, we
specialize it to double structures on rational normal curves. To
every double structure we associate a triple of integers $ (2r,g,n)
$ where $ r $ is the degree of the support, $ n \geq r $ is the
dimension of the projective space containing the double curve, and
$ g $ is the arithmetic genus of the double curve.
We compute also some numerical invariants of the constructed
curves, and we show that the family of double structures with a
given triple $ (2r,g,n) $ is irreducible. Furthermore, we prove
that the general double curve in the families associated to $
(2r,r+1,r) $ and $ (2r,1,2r-1) $ is arithmetically Gorenstein.
Finally, we prove that the closure of the locus containing double
conics of genus $ g \leq -2 $ form an irreducible component of the
corresponding Hilbert scheme, and that the general double conic is
a smooth point of that component. Moreover, we prove that the
general double conic in $ \pp^3 $ of arbitrary genus is a smooth
point of the corresponding Hilbert scheme.
\end{abstract}


\maketitle

\section{Introduction}

Non--reduced projective curves arise naturally when one tries to
classify smooth curves, where a projective curve is a dimension $
1 $ projective scheme without embedded or isolated $
0-$dimensional components. In fact, two of the main tools to
classify projective curves are liaison theory and deformation
theory.

Given two curves $ C $ and $ D $ embedded in the projective space
$ \pp^n,$ we say that they are geometrically linked if they have
no common component and their union is an arithmetically
Gorenstein curve. More than the geometric links, a modern
treatment of the theory takes as its base the algebraic link where
two curves are algebraically linked via the arithmetically
Gorenstein curve $ X $ if $ I_X : I_C = I_D $ and $ I_X : I_D =
I_C,$ where $ I_C, I_D, I_X $ are the saturated ideals that define
the curves $ C, D, X,$ respectively, in the projective space $
\pp^n.$ If $ C $ and $ D $ have no common irreducible component,
the two definitions agree. Liaison theory and even liaison theory
are the study of the equivalence classes of the equivalence
relation generated by the direct link, and by an even number of
direct links, respectively. In $ \pp^3,$ a curve is arithmetically
Gorenstein if, and only if, it is the complete intersection of two
algebraic surfaces. A pioneer in the study of this theory for
curves in $ \pp^3 $ was F. Gaeta (see \cite{g}). In the quoted
paper, he proved that every arithmetically Cohen-Macaulay curve in
$ \pp^3 $ is in the equivalence class of a line. More in general,
every curve sits in an equivalence class, and, for curves in $
\pp^3,$ it is known that every curve in a biliaison class can be
obtained from the curves of minimal degree in the class via a
rather explicit algorithm. This property is known as
Lazarsfeld--Rao property (\cite{mig-book}, Definition 5.4.2), and it was proved in (\cite{mdp}, Ch. IV, Theorem 5.1). The existence of minimal curves and their construction is proved in (\cite{mdp}, Ch. IV, Proposition 4.1, and Theorem 4.3).  In \cite{bbm}, the authors proved the Lazarsfeld--Rao property for the curves in the same biliaison class, without the explicit construction of the minimal curves. The minimal arithmetically Cohen--Macaulay curves are the lines. Hence, the Lazarsfeld--Rao property can be seen as a generalization of Gaeta' s work. Also
if one wants to study smooth curves, the minimal curves in the
biliaison class can have quite bad properties, e.g. they can be
non--reduced, or they can have a large number of irreducible
components. Moreover, the minimal curves in a biliaison class form
an irreducible family of curves with fixed degree and arithmetic
genus. Today, it is not known if the equivalence classes of curves
in $ \pp^n $ have the same properties as those in $ \pp^3 $
(see \cite{nns1}, \cite{cm}, \cite{hart2} for evidence both
ways).

To study the properties of smooth curves, one can also try to
deform the smooth curve to a limit curve and investigate the
properties one is interested in on the limit curve. If those
properties are shared by the limit curve and the deformation
behaves well with respect to the considered properties, then the
general curve shares the same properties of the limit curve.
Often, the limit curves are non--reduced curves. In the papers
\cite{fong}, \cite{bd}, \cite{eg}, the authors study Green's
conjecture concerning the free resolution of a canonical curve by
reducing it to the study of a similar conjecture for double
structures on $ \pp^1 $ called ribbons.

Both described approaches lead to the study of families of curves.
The universal family of curves of fixed degree $ d $ and
arithmetic genus $ g $ is the Hilbert scheme $
\Hilb_{dt+1-g}(\pp^n),$ where, for us, $ \Hilb_{dt+1-g}(\pp^n) $
is the open locus of the full Hilbert scheme corresponding to
locally Cohen-Macaulay $ 1-$dimensional schemes, i.e.
corresponding to curves. Since A. Grothendieck proved its
existence in \cite{gr}, the study of the properties of the Hilbert
scheme attracted many researchers. In spite of their efforts, only a
few properties are known, such as the connectedness of the full
Hilbert scheme proved by R. Hartshorne in \cite{hart1}. A current
trend of research tries to generalize Hartshorne's result on
connectedness to the Hilbert scheme of curves. For partial results
on the problem, see, for example, \cite{hs}, \cite{ns}. In
studying Hilbert schemes, a chance is to relate the local
properties of a point on the Hilbert scheme, e.g. its smoothness
on $ \Hilb_{dt+1-g}(\pp^n),$ and the global properties of the
curve embedded in $ \pp^n.$ With abuse of notation, we denote $ C
$ both the curve in $ \pp^n $ and the corresponding point on $
\Hilb_{dt+1-g}(\pp^n).$ It is well known that the tangent space to
$ \Hilb_{dt+1-g}(\pp^n) $ at a point $ C $ can be identified with
$ H^0(C, \ns_C) $ where $ \ns_C $ is the normal sheaf of the curve
$ C $ as subscheme of $ \pp^n $ (\cite{sernesi}, Theorem 4.3.5).
However, both $ \ns_C $ and its
degree $ 0 $ global sections are far from being well understood
for an arbitrary curve $ C.$

Both in liaison and biliaison theory, and in deformation theory, one has often to consider non--reduced curves. The first general construction for non--reduced curves was given
by D. Ferrand in \cite{fer}, where the author constructs a double
structure on a smooth curve $ C \subset \pp^n.$ The construction
was investigated in \cite{bg}, and generalized in \cite{bafo} to
multiple structures on a smooth support. In the last quoted paper,
the authors present a filtration of a multiple structure $ X $ on
a smooth support $ C $ via multiple structures with smaller
multiplicity on the same curve $ C.$ Moreover, they relate the
properties of $ X $ to the ones of the curves in the filtration. A
different filtration was proposed in \cite{man1}. When the first
multiple structure in either filtration has multiplicity $ 2 $ at
every point, then it comes from Ferrand's construction. In this
sense, double structures are the first step in studying multiple
structures on a smooth curve $ C.$ Because of the previous
discussion, it is interesting to understand if double structures
on curves form irreducible families, if they fill irreducible
components of the Hilbert scheme (if so, they cannot be limit
curves of smooth curves), and if, among them, there are curves
with properties that are preserved under generalization, such as
the property of being arithmetically Gorenstein. In the papers
\cite{nns1}, \cite{nns2}, \cite{nns3}, the authors study the
stated problems for double structures on lines, and more generally
for a multiple structure $ X $ on a line $ L \subset \pp^n $
satisfying the condition $ I_L^2 \subseteq I_X \subseteq I_L,$
called ropes in the literature. In the present paper, we address the
same problems for double structures supported on the most natural
generalization of a line, i.e. a rational normal curve. In
\cite{man2}, the author considers double structures on rational
normal curves, but he is interested in the ones with linear
resolution, a class of curves different from the ones we
investigate.

The plan of the paper is the following. In section
\ref{s-construction}, we recall Ferrand's construction of double
structures on smooth curves, and we specialize it to construct
double structures on rational normal curves. To set notation and
for further use in the paper, we recall some known facts about
rational normal curves. Moreover, we prove that we obtain the
saturated ideal of the double structure directly from the
construction, and we compute the Hilbert polynomial (and hence the arithmetic genus) of the double structure in terms of the numerical data of the construction. Finally, we compute the dimension of the irreducible family of double rational normal curves of given genus.
In section \ref{properties-of-X}, we compute the Hartshorne--Rao function $ h^1
\ii_X(j) $ of such a doubling $ X,$ for $ j \not= 2,$ and we bound
$ h^1 \ii_X(2).$ To get the results, we give also some results
about the ideal sheaf $ \ii_C^2 $ where $ C = X_{red} $ is the
rational normal curve support of $ X.$ Probably, the results we
prove on $ \ii_C^2 $ are folklore, but we did not find references
in literature. In section \ref{agd}, we prove that, among
the double curves we are studying, we can obtain arithmetically
Gorenstein curves. In more detail, it happens in two cases: if $ X
$ has genus $ r+1 $ in $ \pp^{r},$ i.e. $ X $ has degree and genus
of a canonical curve in $ \pp^r,$ and if $ X $ has genus $ 1 $ in
$ \pp^{2r-1},$ i.e. $ X $ has degree and genus of a
non--degenerate normal elliptic curve in $ \pp^{2r-1}.$ The former curves were originally studied in \cite{bd} to understand Green' s conjecture on the free resolution of canonical curves. In the same paper, the authors, together with J. Harris, prove that the considered double structures on rational normal curves are smooth points of the component of the Hilbert scheme containing canonical curves (\cite{bd}, Theorem 6.1). In the
last section of the paper, we study the local properties of $ H(4,
g, n),$ and we show that, if $ g \leq -2,$ then $ H(4, g, n) $ is
open in a generically smooth irreducible component of the Hilbert
scheme $ \Hilb_{4t+1-g}(\pp^n).$ Moreover, we also prove that the
general double conic is a smooth point of $ H(4,g,3) $ with $ g
\geq -1,$ that $ H(4,g,3) $ is not an irreducible component of $
\Hilb_{4t+1-g}(\pp^3),$ and we exhibit the general element of the
irreducible component $ \overline{H(4,g,3)} $ containing $
H(4,g,3).$ The results in this section partially complete the ones
in \cite{ns}. In fact, in \cite{ns}, the authors prove that the
Hilbert scheme $ \Hilb_{4t+1-g}(\pp^3) $ is connected, but do not
study its local properties. By the way, in \cite{ns}, the double
conics are studied as particular curves contained in a double
plane, curves studied in \cite{hs}, and so their construction and
their properties are not considered. Finally, in \cite{bama}, the authors proved that double conics in $ \pp^3_{\mathbb C} $ of genus $ -5 $ are smooth points of the Hilbert scheme $ \Hilb_{4t+6}(\pp^3_{\mathbb C}).$

We want to warmly thank the anonymous referee for his/her comments and remarks and N. Manolache for pointing us a misprint in an earlier draft of the paper.

\section{Construction of double rational normal curves} \label{s-construction}

Let $ K $ be an algebraically closed field of characteristic $ 0 $
and let $ \pp^n $ be the $ n-$dimensional projective space over $
K $ defined as $ \pp^n = \Proj(R := K[x_0, \dots, x_n]).$ If $ X
\subset \pp^n $ is a closed subscheme, we define $ \ii_X $ its
ideal sheaf in $ \oo_{\pp^n} $ and we define the normal sheaf $
\ns_X $ of $ X $ in $ \pp^n $ as $ \ns_X = \homs_{\pp^n}(\ii_X,
\oo_X) = \homs_X(\frac{\ii_X}{\ii_X^2}, \oo_X).$ The saturated
ideal of $ X $ is the ideal $ I_X = \oplus_{j \in \zz} H^0(\pp^n,
\ii_X(j)) \subseteq R,$ and it is a homogeneous ideal. The
homogeneous coordinate ring of $ X $ is defined as $ R_X = R/I_X,$
and it is naturally graded over $ \zz.$ The Hilbert function of $
X $ is then the function defined as $ h_X(j) = \dim_K (R_X)_j,$
degree $ j $ part of $ R_X,$ for $ j \in \zz.$ Finally, it is known that there
exists a polynomial $ P(t) \in \qq[t],$ called Hilbert polynomial
of $ X,$ that verifies $ P(t) = h_X(t) $ for $ t \in \zz, t \gg 0.$
The degree of $ P(t) $ is the dimension of $ X.$ If $ X $ is a
locally Cohen-Macaulay curve, then $ P(t) = dt + 1 - g $ for some
integers $ d, g,$ referred to as degree and arithmetic genus of $
X,$ respectively.

Given a smooth curve $ C \subset \pp^n,$ there is a well known
method, due to D. Ferrand (see \cite{fer}), to construct a
non--reduced curve $ X,$ having $ C $ as support, and multiplicity
$ 2 $ at each point. $ X $ is called a {\em doubling} of $ C.$

Ferrand's method works as follows.

Let $ \ii_C $ be the ideal sheaf of $ C \subset \pp^n,$ and let $
\ii_C/\ii_C^2 $ be its conormal sheaf. If $ \cl $ is a line bundle
on $ C,$ every surjective morphism $ \mu: \frac{\ii_C}{\ii_C^2}
\to \cl $ gives a doubling $ X $ of $ C $ defined by the ideal
sheaf $ \ii_X $ such that $ \ker(\mu) = \ii_X/\ii_C^2.$ The curves
$ C, X $ and the line bundle $ \cl $ are related each other via
the exact sequences \begin{equation} \label{ideal-seq} 0 \to \ii_X
\to \ii_C \to \cl \to 0 \end{equation} and \begin{equation}
\label{structure-seq} 0 \to \cl \to \oo_X \to \oo_C \to
0.\end{equation} Moreover, $ X $ is a locally Cohen-Macaulay curve
and its dualizing sheaf satisfies $ \omega_X \vert C = \cl^{-1}.$

We are interested in studying doublings of rational normal curves,
where, for us, a rational normal curve $ C $ of degree $ r $ is
the image of $$ \pp^1 \stackrel{v_r}\lra \pp^r \lra \pp^n $$ where
$ v_r $ is the Veronese embedding and the second map is a linear
embedding of $ \pp^r $ in $ \pp^n $ with $ r \leq n.$ To make
effective Ferrand's construction in our case, we recall some
known results about rational normal curves, and fix some notation.

Let $ \pp^r \cong L = V(x_{r+1}, \dots, x_n) \subseteq \pp^n =
\Proj(R := K[x_0, \dots, x_n]) $ and let $ C \subset L $ be the
rational normal curve defined by the $ 2 \times 2 $ minors of the
matrix $$ A = \left( \begin{array}{cccc} x_0 & x_1 & \dots &
x_{r-1} \\ x_1 & x_2 & \dots & x_r \end{array} \right).$$ In $ L$
the resolution of the saturated ideal $ I_{C,L} \subset S :=
K[x_0, \dots, x_r] $ of $ C $ is described by the Eagon--Northcott
complex and it is
\begin{equation*} \begin{split} 0 \to & \wedge^r F \otimes
S_{r-2}(G)^* \otimes \wedge^2 G^* \to \wedge^{r-1} F \otimes
S_{r-3}(G)^* \otimes \wedge^2 G^* \to \dots \\ & \dots \to
\wedge^3 F \otimes S_{1}(G)^* \otimes \wedge^2 G^* \to \wedge^2 F
\otimes S_{0}(G)^* \otimes \wedge^2 G^* \to \wedge^0 F \otimes
S_0(G) \to 0 \end{split} \end{equation*} where $ F = S^r(-1), G =
S^2 $ and $ \varphi_A : F \to G $ is defined by the matrix $ A.$

\begin{remark} \label{syz} \rm Because of the definition of $ F $ and $ G $ we
have that the complex ends as follows $$ \dots \to \wedge^3 F
\otimes G^* \stackrel{\varepsilon}{\lra} \wedge^2 F
\stackrel{\phi_A}{\lra} S \to S/I_{C,L} \to 0,$$ where $ \phi_A $
is defined via the $ 2 \times 2 $ minors of $ A.$

Let $ e_1, \dots, e_r $ be the canonical basis
of $ F $ and let $ f_1, f_2 $ be the canonical basis of $ G^*.$
Then, the map $ \varepsilon $ is defined as $$ \varepsilon(e_i
\wedge e_j \wedge e_h \otimes f_k) = x_{i-2+k} e_j \wedge e_h -
x_{j-2+k} e_i \wedge e_h + x_{h-2+k} e_i \wedge e_j $$ for every $
1 \leq i < j < h \leq r, \ k = 1, 2 $ (e.g. see A2.6.1 in \cite{eisenbud-book}).
\end{remark}

Now, we compute the resolution of the saturated ideal $ I_C
\subset R $ of $ C.$ Of course, $ I_C = I_{C,L}^e + I_L =
I_{C,L}^e + \langle x_{r+1}, \dots, x_n \rangle,$ where $ I_{C,L}^e $ is the
extension of $ I_{C,L} $ via the natural inclusion $ S
\hookrightarrow R.$

To get a minimal free resolution of $ I_{C,L}^e $ it suffices to
tensorise by $ \otimes_S R $ the minimal free resolution of $
I_{C,L}.$ To simplify notation, we set $ P_i = \wedge^{i+1} F
\otimes_S S_{i-1}(G)^* \otimes_S \wedge^2 G^* \otimes_S R,$ for $
i = 1, \dots, r-1.$ Hence, the minimal free resolution of $
I_{C,L}^e $ is equal to \begin{equation} \label{res-icle} 0 \to
P_{r-1} \stackrel{\varepsilon_{r-1}}{\lra} P_{r-2} \stackrel{\varepsilon_{r-2}}{\lra} \dots \stackrel{\varepsilon_{3}}{\lra} P_2
\stackrel{\varepsilon_2}{\lra} P_1 \stackrel{\varepsilon_1}{\lra}
I_{C,L}^e \to 0,\end{equation} where the maps are obtained by
tensorising the maps of the minimal free resolution of $ I_{C,L} $
times the identity of $ R.$

The minimal free resolution of $ I_L $ is given by the Koszul
complex over $ x_{r+1}, \dots, x_n.$ Let $ Q = R^{n-r}(-1) $ with
canonical basis $ e_{r+1}, \dots, e_n $ and let $ \delta : Q \to
I_L $ be defined as $ \delta(e_i) = x_i, i = r+1, \dots, n.$ If we
set $ Q_i = \wedge^i Q $ then the minimal free resolution of $ I_L
$ is equal to \begin{equation} \label{res-il} 0 \to Q_{n-r}
\stackrel{\delta_{n-r}}{\lra} Q_{n-r-1} \stackrel{\delta_{n-r-1}}{\lra} \dots \stackrel{\delta_3}{\lra} Q_2
\stackrel{\delta_2}{\lra} Q_1 \stackrel{\delta}{\lra} I_L \to 0
\end{equation} where $ \delta_i = \wedge^i \delta.$

Given the two resolutions (\ref{res-icle}) and (\ref{res-il})
above, we can compute their tensor product (for the definition and
details, see \cite{eisenbud-book}, \S 17.3), and we get
\begin{equation} \label{res-icleil} 0 \to N_{n-1} \lra N_{n-2} \lra
\dots \lra N_2 \lra N_1 \end{equation} where $ N_i = \oplus_{j+k=i+1}
P_j \otimes_R Q_k.$

\begin{lemma} The complex (\ref{res-icleil}) is a minimal
free resolution of $ I_{C,L}^e \cap I_L.$
\end{lemma}

\begin{proof} The free module $ N_i $ is isomorphic to $ R^{\beta_i}(-i-2)
$ with $ \beta_i = \sum_{j+k=i+1} j \binom{r}{j+1} \binom{n-r}k,$
as computed from its definition. Hence, no addendum can be
canceled because of the shifts. It follows that if the complex
(\ref{res-icleil}) is a resolution of $ I_{C,L}^e \cap I_L $ then
it is its minimal free resolution.

At first, we prove that $ I_{C,L}^e \cap \langle x_{r+1}, \dots,
x_l \rangle = I_{C,L}^e \cdot \langle x_{r+1}, \dots, x_l \rangle
$ for every $ l = r+1, \dots, n.$ In fact, if $ f \in I_{C,L}^e
\cap \langle x_{r+1}, \dots, x_l \rangle $ and $ I_{C,L}^e =
\langle g_1, \dots, g_t \rangle,$ then there exist $ h_1, \dots,
h_t \in R $ such that $ f = h_1 g_1 + \dots h_t g_t \in \langle
x_{r+1}, \dots, x_l \rangle.$ For each $ i = 1, \dots, t $ there
exist $ h'_i \in K[x_0, \dots, x_r, x_{l+1}, \dots, x_n] $
and $ h''_i \in \langle x_{r+1}, \dots, x_l \rangle,$ both unique, such that $
h_i = h'_i + h''_i.$ Hence, $ h'_1 g_1 + \dots + h'_t g_t \in
\langle x_{r+1}, \dots, x_l \rangle $ and so it is equal to $ 0,$
because the variables $ x_{r+1}, \dots, x_l $ appear neither in
the $ g_i$'s nor in the $ h'_j$'s. Then, we have $ f = h''_1 g_1 +
\dots + h''_t g_t \in I_{C,L}^e \cdot \langle x_{r+1}, \dots, x_l
\rangle.$

To prove that the complex (\ref{res-icleil}) is a resolution of $
I_{C,L}^e \cap I_L $ we use induction on the number of generators
of $ I_L.$

If $ I_L = \langle x_{r+1} \rangle,$ then $ I^e_{C,L}(-1) \stackrel{\cdot x_{r+1}}{\lra} I^e_{C,L} \cdot I_L $ is a degree $ 0 $ isomorphism, and so the resolution of $ I^e_{C,L} \cdot I_L $ is the one of $ I^e_{C,L} $ shifted by $ -1.$ In this case, $ 0 \to R(-1) \stackrel{\cdot x_{r+1}}{\lra} I_L \to 0 $ is the resolution of $ I_L.$ Thus, the complex (\ref{res-icleil}) is equal to (\ref{res-icle}) tensorised by $ \otimes_R R(-1) $ and so it is the resolution of $ I^e_{C,L} \cdot I_L.$

Assume now that $ I_L = \langle x_{r+1}, \dots, x_n \rangle $ and
that the statement holds for $ I'_L = \langle x_{r+1}, \dots,$ $
x_{n-1} \rangle.$

The two ideals $ I_{C,L}^e \cdot I_L $ and $ I_{C,L}^e \cdot I'_L
+ I_{C,L}^e \cdot \langle x_n \rangle $ are equal. Moreover, $
I_{C,L}^e \cdot I'_L \cap I_{C,L}^e \cdot \langle x_n \rangle $ is
equal to $ I_{C,L}^e \cdot I'_L \cdot \langle x_n \rangle \cong
I_{C,L}^e \cdot I'_L(-1).$ In fact, let $ f \in I_{C,L}^e $ and
assume that $ x_n f \in I_{C,L}^e \cdot I'_l.$ Then, $ x_n f \in
I_{C,L}^e $ and $ x_n f \in I'_L.$ But both the ideals are prime
and $ x_n $ belongs neither to $ I_{C,L}^e $ nor to $ I'_L.$
Hence, $ f \in I_{C,L}^e \cap I'_L = I_{C,L} \cdot I'_L $ and the
statement follows because the converse inclusion is evident. We
have then the following short exact sequence $$ 0 \to I_{C,L}^e
\cdot I'_L(-1) \lra I_{C,L}^e \cdot I'_L \oplus I_{C,L}^e \cdot
\langle x_n \rangle \lra I_{C,L}^e \cdot I_L \to 0 $$ and the
claim follows by applying the mapping cone procedure.
\end{proof}

\begin{proposition} With the same notation as before, the minimal
free resolution of $ I_C $ is $$ 0 \to N_{n-1}
\stackrel{\varepsilon'_{n-1}}{\lra} \dots
\stackrel{\varepsilon'_3}{\lra} P_2 \oplus Q_2 \oplus N_1
\stackrel{\varepsilon'_2}{\lra} P_1 \oplus Q_1
\stackrel{\varepsilon'_1}{\lra} I_C \to 0 $$ where $
\varepsilon'_2 : N_1 = P_1 \otimes Q_1 \lra P_1 \oplus Q_1 $ is
defined as $$ \left( \begin{array}{ccc} -x_{r+1} id_{P_1} & \cdots
& -x_n id_{P_1} \\ \varepsilon_1 & \cdots & \varepsilon_1
\end{array} \right).$$
\end{proposition}

\begin{proof} The ideal $ I_C $ is equal to $ I_{C,L}^e + I_L.$ Hence,
we have the short exact sequence $$ 0 \to I_{C,L}^e \cap I_L \lra
I_{C,L}^e \oplus I_L \lra I_C \to 0.$$ By applying the mapping
cone procedure, we get a free resolution of $ I_C $ that has the
shape of our claim. The minimality of the resolution follows
because a cancelation takes place in the resolution only if a
free addendum of $ N_i $ splits from the map $ N_i \to P_i \oplus
Q_i.$ This cannot happen because $ N_i \cong R^{\beta_i}(-i-2),$ $
P_i \cong R^{i \binom r{i+1}}(-i-1),$ and $ Q_i \cong
R^{\binom{n-r}i}(-i) $ and so the twists do not allow the
splitting of free addenda.
\end{proof}

By sheafifying the previous resolutions, we get the minimal
resolutions of $ \ii_C $ and $ \ii_L $ over $ \oo_{\pp^n},$  and
of $ \ii_{C,L} $ over $ \oo_L,$  that we shall use in what follows. As
standing notation, the map $ \varepsilon'_i $ of the resolution of
$ I_C $ will become $ \tilde\varepsilon'_i $ after sheafifying the
resolution, and the same for the other maps.

As previously explained, to construct a double structure $ X $
supported on $ C $ we need a surjective morphism $ \mu:
\frac{\ii_C}{\ii_C^2} \to \cl $ where $ \cl $ is an invertible
sheaf on $ C.$ We know that $ \ii_C/\ii_C^2 \cong \ii_C \otimes
\oo_C $ and so, if we tensorise the resolution of $ \ii_C $ with $
\oo_C $ we get $ \ii_C/\ii_C^2 \cong \coker(\tilde\varepsilon'_2
\otimes id_{\oo_C}).$ Moreover, it is easy to prove the following
\begin{proposition} $ \coker(\tilde\varepsilon'_2 \otimes id_{\oo_C}) \cong
\coker(\tilde\varepsilon_2 \otimes id_{\oo_C}) \oplus \oo_C^{n-r}(-1).$
\end{proposition}

\begin{proof} The restrictions of $ \tilde\varepsilon'_2 \otimes
id_{\oo_C} $ to $ \tilde N_1 $ and to $ \tilde Q_2 $ are the null
maps because the entries of the mentioned restrictions belong to $
I_C.$
\end{proof}

The curve $ C $ is isomorphic to $ \pp^1,$ and an isomorphism $ j:
\pp^1 \to C $ is defined as $ j(t:u) = (t^r: t^{r-1} u: \dots :
u^r : 0: \dots : 0).$ We have that $ j^*(\oo_C^{n-r}(-1)) =
\oo_{\pp^1}^{n-r}(-r) $ and $ j^*(\coker(\tilde\varepsilon_2
\otimes id_{\oo_C})) = \coker(j^*(\tilde\varepsilon_2 \otimes
id_{\oo_C})) \cong \oo_{\pp^1}^{r-1}(-r-2) $ (see Lemma 5.4 in \cite{bd} for the last isomorphism). Hence, on $ \pp^1,$
the conormal sheaf of $ C \subseteq \pp^n $ is isomorphic to $
\oo_{\pp^1}^{r-1}(-r-2) \oplus \oo_{\pp^1}^{n-r}(-r).$

Now, we make some effort to explicitly write the previous isomorphism.

\begin{lemma} \label{lemma-psi} Let $ e_1, \dots, e_r $ and $ g_1,
\dots, g_{r-1} $ be the canonical bases of $ \oo_{\pp^1}^r(-r) $
and of $ \oo_{\pp^1}^{r-1}(-r-2),$ respectively, and let $ \psi_r:
\wedge^2 \oo_{\pp^1}^r(-r) \lra \oo_{\pp^1}^{r-1}(-r-2) $ be
defined as $$ \psi_r(e_p \wedge e_q) = \sum_{h=p}^{q-1} t^{r-h-1}
u^{h-1} g_{p+q-1-h}, \quad \mbox{ for every } 1 \leq p < q \leq
r.$$ Then, $ \psi_r $ is surjective and $ \ker(\psi_r) \cong \oo_{\pp^1}^{\binom{r-1}2} (-2r-2).$
\end{lemma}

\begin{proof} The map $ \psi_r $ is surjective. In fact, $ \psi_r(e_1 \wedge e_2),
\dots, \psi_r(e_1 \wedge e_r) $ are linearly independent at each
point of $ \pp^1 $ except $ (0:1),$ while $ \psi_r(e_1 \wedge e_r),
\dots, \psi_r(e_{r-1} \wedge e_r) $ are linearly independent at each
point of $ \pp^1 $ except $ (1:0).$ Hence, $ \psi_r $ is surjective
at every point of $ \pp^1 $ and so it is surjective.

Then, we have the following short exact sequence $$ 0 \to \ker(\psi_r) \lra \wedge^2 \oo_{\pp^1}^r(-r) \lra \oo_{\pp^1}^{r-1}(-r-2) \to 0,$$ where $ \ker(\psi_r) $ is a locally free $ \oo_{\pp^1}$-module of rank $ \binom r2 - (r-1) = \binom{r-1}2.$ By (\cite{kempf-book}, Proposition 10.5.1), due to Grothendieck, $$ \ker(\psi_r) \cong \oplus_{i=1}^{\binom{r-1}2} \oo_{\pp^1}(-2r-a_i),$$ for some integers $ 0 \leq a_1 \leq \dots \leq a_{\binom{r-1}2} $ which are uniquely determined by $\ker \psi_r.$

It is an easy check to prove that $$ u^2 e_p \wedge e_q - tu e_p
\wedge e_{q+1} - tu e_{p+1} \wedge e_q + t^2 e_{p+1} \wedge
e_{q+1} \in \ker(\psi_r) $$ for $ p = 1, \dots, r-2 $ and $ q = p+1,
\dots, r-1.$ Of course, if $ q = p+1,$ the third addendum is
missing. Hence, we have $ \binom{r-1}2 $ linearly independent
elements of $ \ker(\psi_r),$ and so there is a subsheaf of $ \ker(\psi_r) $ that is isomorphic to $ \oo_{\pp^1}^{\binom{r-1}2}(-2r-2).$ Hence, the statement holds if we prove that $a_1 = 2$ or equivalently that $\ker \psi_r$ does not contain elements of the form $\sum_{1 \leq p < q \leq r} l_{pq} e_p \wedge e_q$ with $\mathrm{deg}(l_{pq}) \leq 1.$ To see this, we proceed by induction on $r.$

Assume $ r = 3.$ Then, $ \psi_3 $ is represented by the matrix $$ \left( \begin{array}{ccc} t & u & 0 \\ 0 & t & u \end{array} \right).$$ The syzygy module is then generated by $ u^2 e_1 \wedge e_2 - tu e_1 \wedge e_3 + t^2 e_2 \wedge e_3,$ and so the claim holds if $ r = 3.$

Assume the claim holds for $ \psi_{r-1}.$ If $ \sum_{1 \leq p < q \leq r} l_{pq} e_p \wedge e_q \in \ker(\psi_r) $ then, for every $ j = 2, \dots, r,$ the $ l_{pq}$'s satisfy the equation $$ l_{1j} t^{r-2} + l_{1,j+1} t^{r-3}u + \dots + l_{1r} t^{j-2} u^{r-j} + u p_j = 0 $$ where $ p_j \in (l_{23}, \dots, l_{r-1, r}),$ and $ p_2 = 0.$ Then, $ l_{1j} \in (u),$ for every $ j = 2, \dots, r.$ Furthermore, for $ j = 2,$ we have also that $ l_{1r} \in (t).$ By degree reasons, $ l_{1r} = 0.$ By substituting in the equation corresponding to $ j = 2,$ we get that $ t $ can be canceled and so we get that $ l_{1, r-1} \in (t).$ Again by degree reasons, $ l_{1, r-1} = 0.$ By iterating the argument, we get that $ l_{1j} = 0 $ for $ j = 2, \dots, r.$ The map $ \psi_r $ restricted to the span of $ e_2 \wedge e_3, \dots, e_{r-1} \wedge e_r $ is the null map on the first addendum of $ \oo_{\pp^1}^{r-1}(-r-2) $ and $ u \psi_{r-1} $ on the remaining addenda. But then $ \sum_{2 \leq p < q \leq r} l_{pq} e_p \wedge e_q \in \ker(\psi_r) $ implies $ \sum_{2 \leq p < q \leq r} l_{pq} e_p \wedge e_q \in \ker(\psi_{r-1}),$ and we conclude by induction assumption.
\end{proof}

\begin{theorem} \label{map-psi} With the same hypotheses as Lemma \ref{lemma-psi}, the sequence $$ \wedge^3 \oo_{\pp^1}^r(-r) \otimes
(\oo_{\pp^1}^2)^* \stackrel{j^*(\tilde\varepsilon_2 \otimes
id_{\oo_C})}{\lra} \wedge^2 \oo_{\pp^1}^r(-r)
\stackrel{\psi_r}{\lra} \oo_{\pp^1}^{r-1}(-r-2) \to 0 $$ is exact.
\end{theorem}

\begin{proof} We have only to prove that $ \ker(\psi_r) = \mbox{Im}(j^*(\tilde\varepsilon_2 \otimes id_{\oo_C})).$

To start, we verify that $ \psi_r \circ j^*(\tilde\varepsilon_2 \otimes id_{\oo_C}) = 0.$

It is a simple computation and its details are \begin{equation*}
\begin{split} & \psi_r(j^*(\tilde\varepsilon_2 \otimes id_{\oo_C}))(e_i \wedge e_j \wedge e_h \otimes
f_k)) = \\ & = t^{r-i+2-k} u^{i-2+k} \psi_r(e_j \wedge e_h) -
t^{r-j+2-k} u^{j-2+k} \psi_r(e_i \wedge e_h) +  t^{r-h+2-k}
u^{h-2+k} \psi_r(e_i \wedge e_j) = \\ & = t^{r-i+2-k} u^{i-2+k}
(t^{r-h} u^{h-2} g_j + \dots + t^{r-j-1} u^{j-1} g_{h-1}) -
t^{r-j+2-k} u^{j-2+k}(t^{r-h} u^{h-2} g_i + \dots \\ & \dots +
t^{r-i-1} u^{i-1} g_{h-1}) + t^{r-h+2-k} u^{h-2+k}(t^{r-j} u^{j-2}
g_i + \dots + t^{r-i-1} u^{i-1} g_{j-1}) = 0. \end{split}
\end{equation*}

As last step, we must prove that $ \ker(\psi_r) =
\mbox{Im}(j^*(\tilde\varepsilon_2 \otimes id_{\oo_C})).$

In Lemma \ref{lemma-psi}, we proved that $$ u^2 e_p \wedge e_q - tu e_p
\wedge e_{q+1} - tu e_{p+1} \wedge e_q + t^2 e_{p+1} \wedge
e_{q+1} $$ for $ p = 1, \dots, r-2 $ and $ q = p+1,
\dots, r-1,$ generate $ \ker(\psi_r).$

Furthermore, we have the equalities \begin{equation*}
\begin{split} u j^*(\tilde\varepsilon_2 \otimes id_{\oo_C})(e_p \wedge & e_q \wedge e_{q+1}
\otimes f_k) - t j^*(\tilde\varepsilon_2 \otimes id_{\oo_C})(e_{p+1} \wedge e_q \wedge
e_{q+1} \otimes f_k) = \\ = & t^{r-q+1-k} u^{q-2+k} (u^2 e_p
\wedge e_q - tu e_p \wedge e_{q+1} - tu e_{p+1} \wedge e_q + t^2
e_{p+1} \wedge e_{q+1}) \end{split} \end{equation*} for every
admissible $ p < q,$ and so $ (\ker(\psi_r))_P =
(\mbox{Im}(j^*(\tilde\varepsilon_2 \otimes id_{\oo_C}))_P) $ at every point $ P \in \pp^1
\setminus\{(1:0), (0:1)\}.$ At $ (t:u) = (0:1) $ the equality of
the stalks follows from $ j^*(\tilde\varepsilon_2 \otimes id_{\oo_C})(e_p \wedge e_q \wedge
e_r \otimes f_2) = e_p \wedge e_q $ for every $ p = 1, \dots, r-2,
q = p+1, \dots, r-1 $ and the fact that $ (\ker(\psi_r))_{(0:1)} $
is generated by $ e_p \wedge e_q $ with $ p = 1, \dots, r-2, q =
p+1, \dots, r-1.$ Analogously, we get the claim at $ (t:u) = (1:0)
$ by computing $ j^*(\tilde\varepsilon_2 \otimes id_{\oo_C})(e_1 \wedge e_p \wedge e_q
\otimes f_1).$
\end{proof}

When there is no confusion, we' ll write $ \psi $ instead of $ \psi_r.$

Of course, thanks to the isomorphism $ j,$ the map $ \mu:
\frac{\ii_C}{\ii_C^2} \to \cl $ can be written also as $ \mu:
\oo_{\pp^1}^{r-1}(-r-2) \oplus \oo_{\pp^1}^{n-r}(-r) \to
\oo_{\pp^1}(-r-2+a) $ where $ j^*(\cl) = \oo_{\pp^1}(-r-2+a) $ for
some $ a \geq 0,$ because the map $ \mu $ is surjective.

Now, the construction can be rewritten as an algorithm: choose the
map $ \mu,$ and consider the map $$ H^0_*(j_*(\mu \circ (\psi
\oplus id))) : H^0_*(C, \wedge^2 \oo_C^r(-1) \oplus
\oo_C^{n-r}(-1)) \lra H^0_* \cl.$$ Let $ F_1 $ be a free $
H^0_*(C, \oo_C)-$module such that the complex $$ F_1
\stackrel{\nu}{\lra} H^0_*(C, \wedge^2 \oo_C^r(-1) \oplus
\oo_C^{n-r}(-1)) \stackrel{H^0_*(j_*(\mu \circ (\psi \oplus
id)))}{\lra} H^0_* \cl $$ is exact. Let $ N $ be a matrix that
represents the map $ \nu,$ and let $ M $ be a lifting of $ N $
over $ R = H^0_*(\pp^n, \oo_{\pp^n}),$ via the canonical
surjective map $ R \to R/I_C = H^0_*(C, \oo_C).$ The ideal $ I_X $
of the doubling $ X $ is generated by $ I_C^2 + [I_C] M,$ where $
[I_C] $ is a row matrix with entries equal to the generators of $
I_C $ in the same order used to write $ \psi.$

Now, we investigate more deeply the construction. The data we need
to construct such a curve $ X $ are: $(i)$ a rational normal curve
$ C $ of degree $ r $ in its linear span $ L $ embedded in $
\pp^n,$ for some $ n \geq r,$ together with an isomorphism $ j :
\pp^1 \to C;$ $(ii)$ a surjective map $ \mu:
\oo_{\pp^1}^{r-1}(-r-2) \oplus \oo_{\pp^1}^{n-r}(-r) \to
\oo_{\pp^1}(-r-2+a) $ for some $ a \geq 0.$

\begin{remark} \rm For $ r \geq 3 $ and $ a \geq 0
$ there exists always a surjective map $ \mu,$ while, for $ r = n = 2,$ there exists a surjective map
$ \mu $ if, and only if, $ a = 0.$
\end{remark}

\begin{theorem} \label{sat} Let $ X $ and $ X' $ be double structures on two
rational normal curves $ C $ and $ C'.$ Then, $ X = X' $ if, and
only if, $ C = C' $ and the target maps $ \mu $ and $ \mu' $
differ by an automorphism of $ \oo_{\pp^1}(-r-2+a) $ after
changing $ j' $ with $ j.$
\end{theorem}

\begin{proof} Assume first that $ C = C' $ and $ j = j'.$ If $ \mu
$ and $ \mu' $ differ by an automorphism of $ \oo_{\pp^1}(-r-2+a)
$ then the maps $ \mu \circ \psi $ and $ \mu' \circ \psi $ have
the same kernel, and so the curves $ X $ and $ X' $ are defined by
the same ideal, i.e. they are equal each other.

Conversely, if $ X $ and $ X' $ are defined by the same ideal,
then $ X_{red} $ and $ X'_{red} $ are the same curve $ C,$ because
the supporting curve is defined by the only minimal prime ideal
associated to $ I_X.$ Up to compose $ j' $ with an isomorphism of
$ \pp^1 $ we can assume that $ j = j'.$ The claim follows from
(\cite{bafo}, (1.1))
\end{proof}

We show with an example how to compute the ideal of such a
doubling.
\begin{example}\label{ex1} \rm Let $ C \subset \pp^3 = \Proj(K[x,y,z,w]) $ be
the twisted cubic curve whose ideal is $ I_C = (y^2-xz, yz-xw,
z^2-yw).$ We want to construct a double structure $ X $ on $ C $
contained in $ \pp^3.$ As explained, if we set $ a = 1,$ we must
choose a surjective map $$ \mu: \oo_{\pp^1}^2(-5) \to
\oo_{\pp^1}(-4).$$ Set $ \mu = (t,u).$ The map $ \mu \circ \psi $
is given by $ \mu \circ \psi = (t^2, 2tu, u^2),$ while the map $
H^0_*(j_*(\mu \circ \psi)) $ is given by either $ (y, 2z, w) $ or
$ (x, 2y, z) $ (the two apparently different maps agree over $ C
\setminus \{(1:0:0:0), (0:0:0:1) \} $). Of course, to get the two
expressions we multiplied $ \mu \circ \psi $ times $ t $ and $ u $
so that the entries have degree multiple of $ r = 3,$ and then we
used the isomorphism $ j.$ The free $ R/I_C-$module $ F_1 $ that
makes exact the complex $$ F_ 1 \to (R/I_C)^3(-2) \to H^0_* \cl $$
is $ F_1 = (R/I_C)^6(-3) $ and a matrix that represents the map $
F_1 \to (R/I_C)^3(-2) $ is $$ N = \left( \begin{array}{cccccc} 2y
& 2z & 2w & 0 & 0 & 0 \\ -x & -y & -z & y & z & w \\ 0 & 0 & 0 &
-2x & -2y & -2z \end{array} \right).$$ By lifting $ N $ to a
matrix $ M $ over $ R $ via $ R \to R/I_C $ we get $ M = N,$ where
the entries are polynomials in $ R $ and no more equivalence
classes in $ R/I_C.$ The double structure $ X $ is then defined by
$$ I_X = I_C^2 + [I_C] M.$$
\end{example}

\begin{proposition} Let $ C \subseteq \pp^r \cong L \subset \pp^n
$ be a rational normal curve of degree $ r,$ and let $ j: \pp^1
\to C $ be an isomorphism. Let $ \mu: \oo_{\pp^1}^{r-1}(-r-2)
\oplus \oo_{\pp^1}^{n-r}(-r) \to \oo_{\pp^1}(-r-2+a) $ be a
surjective map, and let $ X $ be the double structure on $ C $
associated to $ \mu.$ Then, the ideal $ I_X = I_C^2 + [I_C] M,$
constructed as explained, is saturated.
\end{proposition}

\begin{proof} Let $ J = I_X^{sat}.$ From the inclusion $ J
\subseteq I_C,$ it follows that there exists a matrix $ M' $ such
that $ J = I_C^2 + [I_C] M'.$ Let $ N, N' $ be the images of $ M,
M',$ respectively, when we restrict the last two matrices to $
R/I_C.$ The matrices $ N $ and $ N' $ both present $ H^0_*(C,
\wedge^2 \oo_C^r(-1) \oplus \oo_C^{n-r}(-1))/\ker(H^0_*(j_*(\mu
\circ (\psi \oplus id)))),$ and so the columns of $ N $ (resp. $
N'$) are combination of the ones of $ N' $ (resp. $ N $). Hence, $
I_X = J $ and $ I_X $ is saturated.
\end{proof}

We want to prove some results about families of doublings. Before
stating and proving those results, we compute the Hilbert
polynomial of a doubling $ X $ in terms of the degree of $ \cl.$
Of course, the degree of $ X $ is twice the degree of the rational
normal curve $ C = X_{red} $ and so we have to compute the genus
of $ X.$

\begin{proposition} \label{genus} Let $ X $ be a doubling of a
degree $ r $ rational normal curve $ C $ defined by a map $ \mu $
as above. Then, the Hilbert polynomial of $ X $ is $ P_X(t) = 2rt
+ a -r,$ and so its arithmetic genus $ g_X $ is equal to $ r+1-a.$
\end{proposition}

\begin{proof} By construction, the curves $ C $ and $ X $ and the
invertible sheaf $ \cl $ are related via the short exact sequence
(\ref{structure-seq}), and so the Hilbert polynomial $ P_X(t) $ of
$ X $ is equal to the sum of the Hilbert polynomial $ P_C(t) $ of
$ C $ and of the Euler characteristic $ \chi \cl \otimes
\oo_{\pp^n}(t).$ By restriction to $ \pp^1 $ we get $ \chi \cl
\otimes \oo_{\pp^n}(t) = \chi \oo_{\pp^1}(rt-r-2+a) = rt-r-1+a.$
The Hilbert polynomial of $ C $ is equal to $ P_C(t) = rt+1,$ and
so the claim follows.
\end{proof}

Now, we describe a parameter space for the doublings of the
rational normal curves of fixed degree and genus.

From Proposition \ref{genus}, it follows that if we fix degree and
genus of $ X $ then we fix the degree $ r $ of the rational normal
curve $ C = X_{red} $ and the twist $ a = r+1-g \in \zz $ for the
map $ \mu: \oo_{\pp^1}^{r-1}(-r-2) \oplus \oo_{\pp^1}^{n-r}(-r)
\to \oo_{\pp^1}(-r-2+a) = \oo_{\pp^1}(-1-g).$

Let $ P(t) = 2rt + 1 - g $ be a polynomial and let $
\Hilb_{p(t)}(\pp^n) $ be the Hilbert scheme parameterizing locally
Cohen--Macaulay curves of $ \pp^n $ with Hilbert polynomial $
P(t).$ Let $ H(2r, g, n) $ be the locus in $ \Hilb_{p(t)}(\pp^n) $
whose closed points correspond to double structures of genus $ g $
on smooth rational normal curves of degree $ r $ embedded in $
\pp^n.$ Let $ H(r,n) $ be the locus in $ \Hilb_{rt+1}(\pp^n) $
whose closed points are smooth rational normal curves of degree $
r $ in $ \pp^n.$ $ H(r,n) $ is open in an irreducible component of
$ \Hilb_{rt+1}(\pp^n) $ of dimension $ \dim PGL_r - \dim PGL_1 +
\dim Grass(n-r,n) = (n+1)(r+1) - 4,$ where $ Grass(n-r,n) $ is the
Grassmannian of the linear spaces of dimension $ n-r $ in $
\pp^n.$ Furthermore, there is a natural map $ \varphi: H(2r, g, n)
\to H(r,n) $ defined as $ \varphi(X) = X_{red} $ where, with abuse
of notation, we denote $ X $ both the subscheme in $ \pp^n $ and
the closed point in the Hilbert scheme. The fibers of $ \varphi $
are isomorphic to $$ \Hom(\oo_{\pp^1}^{r-1}(-r-2) \oplus
\oo_{\pp^1}^{n-r}(-r) \to \oo_{\pp^1}(-1-g)) /
\mbox{Aut}(\oo_{\pp^1}(-1-g)) $$ and hence they are irreducible
and smooth of dimension $ (n-1) (r+1-g) + 2r - 2-n.$ Moreover,
both $ H(2r, g, n) $ and $ H(r,n) $ are stable under the action of
$ PGL_n,$ and so we have proved the following
\begin{theorem} \label{hilb-comp} $ H(2r, g, n) $ is irreducible of dimension $$ \dim H(2r,
g, n) = (n+1)(2r+1-g) - 7 + 2g.$$
\end{theorem}

\begin{remark} \rm Of course, we do not know if $ H(2r, g, n) $ is
an irreducible component of $ \Hilb_{P(t)}(\pp^n).$ It is
reasonable that, under suitable hypotheses on $ g,$ it is so. In
the last section of the paper, we will study the local properties
of $ H(4, g, n) $ for $ n \geq 3.$
\end{remark}

\begin{corollary} \label{dim-comp} Every irreducible
component of $ \Hilb_{2rt+1-g}(\pp^n) $ containing $ H(2r,g,n) $ has dimension
$ \geq (n+1)(2r+1-g) - 7 + 2g.$
\end{corollary}


\section{Cohomology estimates} \label{properties-of-X}

In this section, we show how to compute the Rao function of a
double rational normal curve in terms of the data of the
construction. To achieve the results, we have to investigate also
the curve $ D $ defined by the ideal sheaf $ \ii_C^2.$

At first, we restate a deep result about $ D $ due to
J. Wahl (see \cite{wahl}, Theorem 2.1), that holds more generally
for every Veronese embedding of a projective space.
\begin{theorem} Let $ C $ be a rational normal curve. Then, $ H^1 \ii_C^2(j)
= 0 $ for $ j \not= 2.$
\end{theorem}

We want to compute $ h^1 \ii_C^2(2) $ and the ideal $ I_D =
(I_C^2)^{sat}.$ The following results are probably known in
literature, but we add their proofs for completeness. The
computation of the generators of $ I_D $ rests on some direct
calculations and on some basic results on initial ideals. Let us
recall the following
\begin{lemma} \label{in-id}
Let $ I, J $ be ideals in a polynomial ring $ R.$ If $ I \subseteq
J $ and $ \mbox{in}(I) \supseteq \mbox{in}(J),$ then $ I = J,$ no
matter what term ordering we use to compute the initial ideal.
\end{lemma}

\begin{proof} If $ I \subseteq J $ then $ \mbox{in}(I) \subseteq
\mbox{in}(J).$ From our hypothesis, it follows that $ \mbox{in}(I)
= \mbox{in}(J),$ and the claim follows from (\cite{eisenbud-book},
Lemma 15.5).
\end{proof}

Now, we prove some results about the initial ideal of the ideal of
a rational normal curve and of its square.
\begin{lemma} \label{hf-in} Let $ C \subset \pp^r $ be the rational normal
curve of degree $ r $ generated by the $ 2 \times 2 $ minors of
the matrix $$ A = \left( \begin{array}{cccc} x_0 & x_1 & \dots &
x_{r-1} \\ x_1 & x_2 & \dots & x_r \end{array} \right).$$ Then,
with respect to the degrevlex ordering of the terms in $ R,$ we
have that \begin{enumerate} \item $ \mbox{in}(I_C) = \langle x_1, \dots,
x_{r-1} \rangle^2;$ \item $ \mbox{in}(I_C^2) = \langle x_1, \dots, x_{r-1} \rangle^4 +
x_0 \langle x_2, \dots, x_{r-1} \rangle^3 + x_r \langle x_2, \dots, x_{r-2} \rangle^3;$ \item
$ (\mbox{in}(I_C^2))^{sat} = \mbox{in}(I_C^2) + \langle x_2, \dots,
x_{r-2} \rangle^3.$
\end{enumerate} \end{lemma}

\begin{proof} Let $ f_{ij} $ be the minor given by the $
i-$th and $ j-$th columns of $ A,$ i.e. $ f_{ij} = x_{i-1} x_j -
x_i x_{j-1},$ with $ 1 \leq i < j \leq r.$ The leading term of $
f_{ij},$ with respect to the graded reverse lexicographic order, is $
\mbox{in}(f_{ij}) = x_i x_{j-1},$ and so $ \mbox{in}(I_C)
\supseteq \langle x_1, \dots, x_{r-1} \rangle^2.$ On the other hand, it is easy
to verify that $ \{f_{12}, \dots, f_{r-1,r}\} $ is a Gr\"obner
basis of $ I_C $ and so the first claim holds.

Let $ D \subset \pp^r $ be the curve defined by the ideal $ I_D =
H^0_*(\ii_C^2) = (I_C^2)^{sat}.$ Thanks to the previous Theorem,
the homogeneous parts of the ideals $ I_C $ and $ I_D $ are
related each other from the exact sequence $$ 0 \to (I_D)_t \to
(I_C)_t \to H^0(\oo_{\pp^1}^{r-1}(-r-2+rt)) \to 0 $$ for every $ t
\geq 3.$ Hence, $ \dim_K (I_D)_t = \binom{r+t}r -(r^2 t + 2 -
r^2).$ It is well known that $ \dim_K (I)_t = \dim_K
(\mbox{in}(I))_t $ for every homogeneous ideal
(\cite{eisenbud-book}, Theorem 15.26). So, $ \dim_K
(\mbox{in}(I_C^2))_t \leq \binom{r+t}r -(r^2 t + 2 - r^2).$

Consider the monomial ideal $ J = \langle x_1, \dots, x_{r-1} \rangle^4 + x_0
\langle x_2, \dots, x_{r-1} \rangle^3 + x_r \langle x_2, \dots, x_{r-2} \rangle^3.$

We claim that $ J = \mbox{in}(I_C^2).$

It is a straightforward computation to check that $$ \dim_K (J)_t
= \binom{r+t}r - (r^2 t + 2 - r^2),$$ for $ t \geq 4.$ For
example, it is easy to enumerate the degree $ t $ monomials not in
$ J.$ It follows that $ \dim_K (\mbox{in}(I_C^2))_t \leq \dim_K
(J)_t.$ So, if $ \mbox{in}(I_C^2) \supseteq J,$ then $
\mbox{in}(I_C^2) = J.$ Furthermore, we get also that $ (I_D)_t =
(I_C^2)_t $ for $ t \geq 4,$ because their homogeneous parts have
the same dimension for $ t \geq 4.$

Of course, $ \mbox{in}(I_C^2) \supset \langle x_1, \dots, x_{r-1} \rangle^4,$
because of the first claim. Moreover, if $ 2 \leq i \leq j \leq h
\leq r-1,$ the leading term of $ f_{1i} f_{j,h+1} - f_{1,j+1}
f_{i-1,h+1} $ is equal to $ x_0 x_i x_j x_h $ and so $ x_0 \langle x_2,
\dots, x_{r-1} \rangle^3 \subset \mbox{in}(I_C^2).$ Finally, if $ 2 \leq
i \leq j \leq h \leq r-2,$ then the leading term of $ f_{i,j+1}
f_{h+1,r} - f_{i,h+2} f_{j,r} + f_{i-1,r} f_{j+1, h+2} $ is equal
to $ x_i x_j x_h x_r $ and so $ x_r \langle x_2, \dots, x_{r-2} \rangle^3
\subset \mbox{in}(I_C^2).$ Because of the previous argument, the
second claim follows. The equality in $ (3) $ follows from the
definition of saturation.
\end{proof}

Now, we can compute the generators of $ I_D.$
\begin{theorem} With the same hypotheses as before, let $ I' $ be
the ideal generated by the $ 3 \times 3 $ minors of the matrix $$
B = \left( \begin{array}{cccc} x_0 & x_1 & \dots & x_{r-2} \\ x_1
& x_2 & \dots & x_{r-1} \\ x_2 & x_3 & \dots & x_r \end{array}
\right).$$ Then, $ I_D = I_C^2 + I'.$
\end{theorem}

\begin{proof} By construction, $ \mbox{in}(I') \supseteq \langle x_2,
\dots, x_{r-2} \rangle^3,$ and by Lemma \ref{hf-in}(3), we have $
(\mbox{in}(I_C^2))^{sat} $ $ = \mbox{in}(I_C^2) + \langle x_2, \dots,
x_{r-2} \rangle^3.$ Moreover, we have the following chain of inclusions
$$ \mbox{in}(I_C^2 + I') \supseteq \mbox{in}(I_C^2) +
\mbox{in}(I') \supseteq \mbox{in}(I_C^2) + \langle x_2, \dots, x_{r-2} \rangle^3
= (\mbox{in}(I_C^2))^{sat} \supseteq \mbox{in}((I_C^2)^{sat}) =
\mbox{in}(I_D).$$ So, from the inclusion $ I_C^2 \subseteq I_D $
and Lemma \ref{in-id}, it is enough to show that $ I' \subset
I_D.$

Let us consider the integers $ i,j,h $ with $ 2 \leq i < j < h
\leq r $ and let $$ g_{ijh} = \det \left( \begin{array}{ccc}
x_{i-2} & x_{j-2} & x_{h-2} \\ x_{i-1} & x_{j-1} & x_{h-1} \\ x_i
& x_j & x_h \end{array} \right).$$ It is evident that the two following
equalities hold $$ g_{ijh} = x_{i-2} f_{jh} - x_{j-2} f_{ih} +
x_{h-2} f_{ij} = x_i f_{j-1,h-1} - x_j f_{i-1,h-1} + x_h
f_{i-1,j-1}.$$

We want to prove that $ x_k g_{ijh} \in I_C^2 $ for every $ k = 0,
\dots, r.$

From Proposition \ref{syz}, we know that
$$ x_{i-1} f_{jh} - x_{j-1} f_{ih} + x_{k-1} f_{ij} = x_i f_{jh} -
x_j f_{ih} + x_h f_{ij} = 0.$$ Hence, the following easy
computations prove the claim $$ x_k g_{ijh} = \left\{
\begin{array}{cl} f_{i-1,k} f_{jh} - f_{j-1,k} f_{ih} + f_{h-1,k}
f_{ij} & \mbox{if } k \geq h \\ f_{i,k+1} f_{j-1,k-1} - f_{j,k+1}
f_{i-1,h-1} + f_{k+1,h} f_{i-1,j-1} & \mbox{if } k < h \end{array}
\right. $$ where, in the last equation, we use the convention that
$ f_{aa} = 0 $ for every $ a,$ and $ f_{ab} = f_{ba} $ if $ a > b.$
\end{proof}

\begin{remark} \rm If $ C $ is a line or a smooth conic in $ \pp^n
$ then $ I_C^2 $ is generated by $ \binom r2 $ polynomials. By the
way, those two cases are the only for which $ I_C $ is a complete
intersection ideal. If $ C $ is a twisted cubic curve, then $
I_C^2 $ is saturated. If $ r \geq 4,$ then $ I_C^2 $ is no more
saturated, but $ I_C^2 $ and its saturation agree from degree $ 4
$ on.
\end{remark}

\begin{remark} \rm In the proof of Theorem \ref{map-psi}, we
checked that $ \psi(j^*(\tilde\varepsilon_2 \otimes
id_{\oo_C})(e_i \wedge e_j \wedge e_h \otimes f_k)) = 0,$ and so $
j^*(\tilde\varepsilon_2 \otimes id_{\oo_C})(e_i \wedge e_j \wedge
e_h \otimes f_k) \in \ker( \mu_\vert \circ \psi) $ where $
\mu_{\vert} $ is the restriction of $ \mu $ to $ \wedge^2
\oo_{\pp^1}^{r-1}(-r-2).$ The element $ j^*(\tilde\varepsilon_2
\otimes id_{\oo_C})(e_i \wedge e_j \wedge e_h \otimes f_k) $
corresponds to the minor $ g_{ijh} $ and so $ (I_C^2)^{sat}
\subseteq I_X $ for every double structure $ X.$
\end{remark}

Now, we consider the case $ C \subset \pp^r \subset \pp^n.$ As in
section \ref{s-construction}, we denote $ S = K[x_0, \dots, x_r],$
and $ R = K[x_0, \dots, x_n].$
\begin{proposition} Let $ C \subset L \cong \pp^r \subset \pp^n $ be a
rational normal curve. Let $ I_{C,L}, I_C $ be the ideals of $ C $
as a subscheme of $ L $ and of $ \pp^n,$ respectively, and let $
I_L $ be the ideal of $ L $ in $ \pp^n.$  Then, $ (I_C^2)^{sat} =
((I_{C,L}^2)^{sat})^e + (I_{C,L}^e) \cdot I_L + I_L^2 =
((I_{C,L}^2)^{sat})^e + I_C \cdot I_L$ where the extension is via
the natural inclusion $ S \hookrightarrow R.$
\end{proposition}

\begin{proof} The ideals $ I_L^2 $ and $ (I_{C,L}^e) \cdot I_L $ are
obviously contained in $ I_C^2 $ and hence in $ (I_C^2)^{sat}.$
Furthermore, to check the inclusion $ ((I_{C,L}^2)^{sat})^e
\subset (I_C^2)^{sat} $ it is enough to verify that $ f \in
(I_C^2)^{sat} $ for every $ f \in ((I_{C,L}^2)^{sat})^e \cap S.$
Let $ f $ be a homogeneous polynomial in $ ((I_{C,L}^2)^{sat})^e
\cap S.$ It is easy to check that, for every $ i = 0, \dots, n,$
there exists $ m_i \in \nn $ such that $ x_i^{m_i} f \in (I_C^2).$
In fact, if $ 0 \leq i \leq r,$ then there exists $ m_i \in \nn $
such that $ x_i^{m_i} f \in (I_{C,L}^2)^e \cap S \subseteq
(I_C^2),$ while, for $ r+1 \leq i \leq n,$ $ x_i^2 f \in I_C^2 $
because $ I_L^2 \subset I_C^2.$ Hence, $ ((I_{C,L}^2)^{sat})^e
\subseteq (I_C^2)^{sat} $ and the inclusion $
((I_{C,L}^2)^{sat})^e + (I_{C,L}^e) \cdot I_L + I_L^2 \subseteq
(I_C^2)^{sat} $ follows.

To prove the inverse inclusion, let $ f \in (I_C^2)^{sat} $ be a
homogeneous polynomial. There exist $ f_1 \in S $ and $ f_2 \in
I_L $ such that $ f = f_1 + f_2 $ and the decomposition is unique.
For every $ i = 0, \dots, n,$ there exists $ m_i \in \nn $ such
that $ x_i^{m_i} f \in I_C^2.$ Assume $ 0 \leq i \leq r.$ We know
that $ I_C = I_{C,L}^e + I_L,$ and so $ I_C^2 = (I_{C,L}^2)^e +
I_L \cdot I_C.$ Hence, there exist $ g_1 \in (I_{C,L}^2)^e \cap S
$ and $ g_2 \in I_L \cdot I_C $ such that $ x_i^{m_i} f =
x_i^{m_i} f_1 + x_i^{m_i} f_2 = g_1 + g_2.$ It is evident that $
x_i^{m_i} f_1 - g_1 = g_2 - x_i^{m_i} f_2 \in S \cap I_L = 0,$ and
so $ x_i^{m_i} f_1 \in (I_{C,L}^2)^e $ and $ x_i^{m_i} f_2 \in I_L
\cdot I_C.$ By definition of saturation, $ f_1 \in
((I_{C,L}^2)^e)^{sat}.$ By assumption, $ f_2 \in I_L,$ and so $
f_2 = x_{r+1} f_{2,r+1} + \dots + x_n f_{2,n}.$ Then, $ x_i^{m_i}
f_{2,j} \in I_C $ for every $ j = r+1, \dots, n.$ The ideal $ I_C
$ is a prime ideal and $ x_i^{m_i} \notin I_C.$ Hence, $ f_{2,j}
\in I_C,$ $ f_2 \in I_L \cdot I_C $ and the proof is complete.
\end{proof}

A consequence of the previous results is that we can compute also
$ h^1 \ii_C^2(2).$
\begin{corollary} With the same hypotheses as before, $ h^1
\ii_C^2(2) = \binom{r-1}2.$
\end{corollary}

\begin{proof} If we tensor the exact sequence $$ 0 \to \ii_C^2 \to
\ii_C \to \oo_{\pp^1}^{r-1}(-r-2) \oplus \oo_{\pp^1}^{n-r}(-r) \to
0 $$ by $ \oo_{\pp^n}(2) $ and we take the cohomology, we get $$ h^1
\ii_C^2(2) = h^0 \ii_C^2(2) - h^0 \ii_C(2) + (r-1) h^0
\oo_{\pp^1}(r-2) + (n-r) h^0 \oo_{\pp^1}(r) = \binom{r-1}2.$$
\end{proof}

Thanks to the results on $ D,$ we can compute, or at least bound,
the Rao function $ h^1 \ii_X(j), j \in \zz,$ of $ X $ in terms of
the map $ \mu $ which describes the schematic structure of $ X.$

\begin{proposition}\label{rao-f} With the same notation as above,
let $ I_{\mu} $ be the ideal generated by the entries of $ \mu.$
Then, it holds \begin{equation} \label{rao-f-formula} h^1 \ii_X
(j) = \dim_K \left( \frac{K[t,u]}{I_{\mu}} \right)_{rj-r-2+a}
\end{equation} for every $ j \not= 2.$ Moreover, $$ h^1 \ii_X(2)
\leq \dim_K \left( \frac{K[t,u]}{I_{\mu}} \right)_{r-2+a} +
\binom{r-1}2.$$
\end{proposition}

\begin{proof} By construction, we have the short exact sequences $$
0 \to \ii_X \to \ii_C \to \cl \to 0 ,$$ and $$ 0 \to
\frac{\ii_X}{\ii_C^2} \to \frac{\ii_C}{\ii_C^2} \to \cl \to 0
$$ where the first map of them both is the inclusion.

The two sequences fit into the larger commutative diagram
$$ \begin{array}{ccccccc} & 0 & & 0 & & & \\ & \downarrow & &
\downarrow & & & \\ & \ii_C^2 & = & \ii_C^2 & & & \\ & \downarrow
& & \downarrow & & & \\ 0 \to & \ii_X & \lra & \ii_C & \lra & \cl
& \to 0 \\ & \downarrow & & \downarrow & & \parallel & \\ 0 \to &
\frac{\ii_X}{\ii_C^2} & \lra & \frac{\ii_C}{\ii_C^2} & \lra & \cl
& \to 0. \\ & \downarrow & & \downarrow & & & \\ & 0 & & 0 & & &
\end{array} $$ If we twist by $ \oo_{\pp^n}(j) $ and take the
cohomology, we get $$ \begin{array}{ccccccccc} & 0 & & 0 & & & & & \\
& \downarrow & & \downarrow & & & & & \\ & H^0 \ii_C^2(j) & = &
H^0 \ii_C^2(j) & & & & &
\\ & \downarrow & & \downarrow & & & & & \\ 0 \to & H^0 \ii_X(j) & \lra &
H^0 \ii_C(j) & \lra & H^0 \cl(rj) & \lra & H^1 \ii_X (j) & \to 0 \\
& \downarrow & & \downarrow & & \parallel & & \downarrow & \\ 0 \to & H^0
\frac{\ii_X}{\ii_C^2}(rj) & \lra & H^0 \frac{\ii_C}{\ii_C^2}(rj) &
\lra & H^0 \cl(rj) & \lra & \coker_j & \to 0. \\ & \downarrow & & \downarrow & & & \\
& H^1 \ii_C^2(j) & \lra & H^1 \ii_C^2(j) & & & & & \\ & & &
\downarrow & & & & \\ & & & 0 & & & & \end{array} $$

If $ j \not= 2,$ then $ H^1 \ii_C^2(j) = 0 $ and so $$ \coker(H^0
\ii_X(j) \to H^0 \ii_C(j)) = \coker(H^0 \frac{\ii_X}{\ii_C^2}(rj)
\to H^0 \frac{\ii_C}{\ii_C^2}(rj)) $$ as subspaces of $ H^0
\cl(rj).$ Hence, $ H^1 \ii_X(j) \cong \coker_j \cong \left(
\frac{K[t,u]}{I_{\mu}} \right)_{rj-r-2+a}.$

If $ j = 2,$ we set $ \mathcal{A} = \ker(H^0 \cl(2r) \to H^1
\ii_X(2)) $ and $ \mathcal{B} = \ker(H^0 \cl(2r) \to \coker_2).$
Then, the identity of $ H^0 \cl(2r) $ induces an injective map $
\mathcal{A} \to \mathcal{B},$ and the diagram $$
\begin{array}{ccccccc} 0 \to & H^0 \ii_X(2) & \lra & H^0 \ii_C(2)
& \lra & \mathcal{A} & \to 0 \\ & \downarrow & & \downarrow & &
\downarrow & \\ 0 \to & H^0 \frac{\ii_X}{\ii_C^2}(2r) & \lra & H^0
\frac{\ii_C}{\ii_C^2}(2r) & \lra & \mathcal{B} & \to 0
\end{array} $$ induces a surjective map $$ H^1 \ii_C^2(2) \to
\coker(\mathcal{A} \to \mathcal{B}) \cong \ker(H^1 \ii_X(2) \to
\coker_2).$$ Hence, the claim follows from the surjectivity of the
map $ H^1 \ii_X(2) \to \coker_2.$
\end{proof}

\begin{example} [Example \ref{ex1} revisited] \rm The genus of the
curve $ X $ is $ g_X = 4 - 1 = 3,$ because $ a = 1 $ (see Proposition \ref{genus}).
The map $ \mu : \oo_{\pp^1}^2(-5) \to \oo_{\pp^1}(-4) $ was
defined as $ \mu=(t,u) $ and so the Hilbert function of its
cokernel is $$ \dim_K \left( \frac{K[t,u]}{(t,u)} \right)_h =
\left\{
\begin{array}{cl} 1 & \mbox{if } h = 0 \\ 0 & \mbox{otherwise}
\end{array} \right.$$ By Proposition \ref{rao-f}, the Rao function
of $ X $ is equal to
$$ h^1 \ii_X(j) = \dim_K \left( \frac{K[t,u]}{(t,u)} \right)_{3j-4} = 0 $$ for every $ j
\not= 2.$ To compute $ h^1 \ii_X(2) $ we consider the hyperplane $
H = V(w) $ that is general for $ X,$ and the exact sequence $$ 0
\to \ii_X(-1) \to \ii_X \to \ii_{X \cap H \vert H} \to 0.$$ If we
tensorize by $ \oo_{\pp^3}(2) $ and take the cohomology, we get $$ 0
\to H^1 \ii_X(2) \to H^1 \ii_{X \cap H \vert H}(2) \to 0,$$
because $ H^2 \ii_X(1) = 0.$ It is easy to verify that $ h^1
\ii_{X \cap H \vert H}(2) = h^0 \ii_{X \cap H \vert H}(2) $ and so
$ h^1 \ii_X(2) = 0 $ if, and only if, $ X \cap H $ is not
contained in any conic of $ H.$ But $ \left( \frac{I_X + wR}{wR}
\right)^{sat} = \langle 2y^3 - 3 xyz, y^2z-2xz^2, yz^2, z^3
\rangle $ and so $ X $ is an arithmetically Cohen-Macaulay curve
in $ \pp^3.$
\end{example}


\section{Arithmetically Gorenstein double rational normal
curves} \label{agd}

In this section, we want to describe the arithmetically Gorenstein
curves among the double structures on rational normal curves. At
first, we characterize the possible triples $ (2r,g,n) $ and then
we study the possible cases one at a time.

To start, we recall the definition of arithmetically Gorenstein curve.
\begin{definition} A curve $ X \subset \pp^n $ is arithmetically Gorenstein if
its homogeneous coordinate ring $ R_X $ is a Gorenstein ring, or,
equivalently, if $ R_X $ is Cohen-Macaulay and its canonical sheaf
$ \omega_X $ is a twist of the structure sheaf.
\end{definition}

Now, we look for triples $ (2r,g,n) $ for which the property of
being arithmetically Gorenstein is allowed.

\begin{proposition} Let $ C \subset L \cong \pp^r \subset \pp^n, n \geq 3,$
be a rational normal curve of degree $ r,$ let $ \mu :
\oo^{r-1}_{\pp^1}(-r-2) \oplus \oo^{n-r}_{\pp^1}(-r) \to
\oo_{\pp^1}(-1-g) $ be a surjective map, and let $ X $ be the
double structure on $ C $ defined by $ \mu.$ If $ X $ is a
non-degenerate arithmetically Gorenstein curve, then either $
(2r,g,n) = (2r,r+1,r) $ or $ (2r,g,n) = (2r,1,2r-1).$
\end{proposition}

\begin{proof} If $ X $ is an arithmetically Gorenstein curve, then
the second difference $ \Delta^2 h_X $ of its Hilbert function $
h_X $ is a symmetric function. Moreover, if $ X $ is
non-degenerate, then $ \Delta^2 h_X(1) = n-1.$ We have the
equality $ 2r = \deg(X) = \sum_{j=0}^{\infty} \Delta^2 h_X(j) $
and so we get $ r \leq n \leq 2r,$ where the first inequality
comes from the general setting, and the second one from $ \Delta^2
h_X(0) = 1.$ If $ X $ is an arithmetically Gorenstein curve, then
$ \Delta^2 h_X $ is the Hilbert function of the Artinian ring $
R/\langle I_X, h_1, h_2 \rangle $ where $ h_1, h_2 $ are two
linear forms, general with respect to $ X.$ In particular, if $
\Delta^2 h_X (j) = 0 $ for some $ j > 0,$ then $ \Delta^2 h_X (k)
= 0 $ for every $ k \geq j.$ From the above discussion and
inequalities, we get that there are either $ 3 $ or $ 4 $
non--zero entries in $ \Delta^2 h_X.$ In the first case, then $
\Delta^2 h_X = (1, 2r-2, 1) $ and $ (2r, g, n) = (2r, 1, 2r-1).$
In the second case, then $ \Delta^2 h_X = (1, r-1, r-1, 1) $ and $
(2r, g, n) = (2r, r+1, r).$
\end{proof}

\begin{remark} \rm If $ X $ is a double conic in $ \pp^2,$ then $
I_X = \langle q^2 \rangle $ where $ q = xz-y^2 $ defines the
smooth conic that supports $ X.$ $ X $ is arithmetically
Gorenstein with $ \omega_X = \oo_X(1).$
\end{remark}

Now, we characterize the arithmetically Gorenstein double curves
among the ones we can construct with given triple $ (2r,g,n).$

\begin{theorem} \label{canonical} Let $ (2r,g,n) = (2r,r+1,r).$ For every non--zero map $
\mu: \oo^{r-1}_{\pp^1}(-r-2) \to \oo_{\pp^1}(-r-2) $ we get a
non--degenerate arithmetically Cohen--Macaulay curve $ X.$ Furthermore, if $ \mu $ is general, then $ X $ is arithmetically Gorenstein. If $ \mu = (\alpha^{r-2}, $ $ \alpha^{r-3} \beta, \dots, \beta^{r-2}) $ then
$ X $ is contained in a cone over a rational normal curve $ C' \subset \pp^{r-1}.$
\end{theorem}

\begin{proof} Let $ \mu = (a_0, \dots, a_{r-2}) \not= 0 $ with $ a_i \in K.
$ Let $ X $ be the curve we get by doubling a rational normal
curve $ C $ via $ \mu.$

We prove that $ X $ is an arithmetically Cohen-Macaulay
curve. By Proposition \ref{rao-f}, the surjectivity of $ \mu $
implies that $ h^1 \ii_X(j) = 0 $ for $ j \not=2,$ and so we have
to prove that $ h^1 \ii_X(2) = 0,$ too. The map $ \mu \circ \psi $
can be written also as the composition $ \mu' \circ \psi' $ where
$ \mu' : \oo_{\pp^1}^{r-1}(-2r) \to \oo_{\pp^1}(-r-2) $ is defined
as $ \mu' = (t^{r-2}, t^{r-3} u, \dots, u^{r-2}) $ and $ \psi' :
\wedge^2 \oo_{\pp^1}^r(-r) \to \oo_{\pp^1}^{r-1}(-2r) $ is defined
by the matrix obtained from the one of $ \psi $ by substituting $
t^i u^j, i+j = r-2,$ with $ a_j,$ where $ \psi $ was defined in
Theorem \ref{map-psi}. The matrix of $ \psi' $ is full rank for
whatever non--zero map $ \mu,$ and so it has exactly $ \binom r2 -
(r-1) = \binom{r-1}2 $ linearly independent degree $ 0 $ syzygies.
Furthermore, a degree $ 0 $ syzygy of $ \mu \circ \psi $ is a
degree $ 0 $ syzygy of $ \psi' $ and hence, $ h^0 \ii_X(2) =
\binom{r-1}2.$ From the exact sequence
$$ 0 \to H^0 \ii_X(2) \to H^0 \ii_C(2) \to H^0
\frac{\ii_C}{\ii_X}(2) \to H^1 \ii_X(2) \to 0,$$ and from $ h^0
\ii_C(2) = \binom r2, h^0 \frac{\ii_C}{\ii_X}(2) = h^0
\oo_{\pp^1}(r-2) = r-1, h^0 \ii_X(2) = \binom{r-1}2 $ we get that
$ h^1 \ii_X(2) = 0.$

The curve $ C $ is rational normal and so $ \mbox{Pic}(C) = \zz.$ We know that
the line bundle $ \cl $ verifies $ j^*(\cl) = \oo_{\pp^1}(-r-2) $
and so $ \cl = \omega_C(-1) $ where $ \omega_C $ is the canonical
sheaf of $ C.$ Hence, the curve $ X $ is defined via the exact
sequence $$ 0 \to \ii_X \to \ii_C \to \omega_C(-1) \to 0.$$ From the exact sequence, we can compute the Hilbert function of $ X $ and its second difference. In particular, we get $ \Delta^2 h_X = (1, r-1, r-1, 1),$ as expected.

Assume now that $ \mu = (\alpha^{r-2}, \alpha^{r-3} \beta, \dots,
\beta^{r-2}) $ for some $ (\alpha, \beta) \in K^2 \setminus \{
(0,0) \}.$ The $ 2 \times 2 $ minors of the matrix $$ \left(
\begin{array}{cccc} \beta x_0 - \alpha x_1 & \beta x_1 - \alpha
x_2 & \dots & \beta x_{r-2} - \alpha x_{r-1} \\ \beta x_1 - \alpha
x_2 & \beta x_2 - \alpha x_3 & \dots & \beta x_{r-1} - \alpha
x_{r} \end{array} \right) $$ define a cone in $ \pp^r $ over a
rational normal curve of $ \pp^{r-1}.$ We want to prove that they
belong to $ I_X.$ To this end, let $ 1 \leq i < j \leq r-1.$ The
minor $ F_{i,j} $ given by the $ i-$th and $ j-$th columns is
equal to $$ F_{i,j} = \det \left( \begin{array}{cc} \beta x_{i-1}
- \alpha x_i & \beta x_{j-1} - \alpha x_j \\ \beta x_i - \alpha
x_{i+1} & \beta x_j - \alpha x_{j+1} \end{array} \right) = \beta^2
f_{ij} - \alpha \beta f_{i,j+1} - \alpha \beta f_{i+1,j} +
\alpha^2 f_{i+1,j+1} $$ where $ f_{pq} = x_{p-1} x_q - x_p x_{q-1}
$ is a generator of $ I_C.$ The claim follows if we prove that $$
\beta^2 e_i \wedge e_j - \alpha \beta e_i \wedge e_{j+1} - \alpha
\beta e_{i+1} \wedge e_j + \alpha^2 e_{i+1} \wedge e_{j+1} $$ is a
syzygy of $ \mu \circ \psi $ for every $ 1 \leq i < j \leq r-2.$
Since $ \mu(g_h) = \alpha^{r-h-1} \beta^{h-1}, h = 1, \dots, r-1,$
where $ g_1, \dots, g_{r-1} $ is the canonical basis of $
\oo_{\pp^1}^{r-1}(-r-2),$ then, we have \begin{equation*}
\begin{split} \mu \circ \psi (\beta^2 e_i \wedge e_j - \alpha \beta
& e_i \wedge e_{j+1} - \alpha \beta e_{i+1} \wedge e_j + \alpha^2
e_{i+1} \wedge e_{j+1}) =\\ = \mu(t^{r-j-1} u^{j-1} &(-\alpha
\beta g_i + \alpha^2 g_{i+1}) + t^{r-j} u^{j-2} (\beta^2 g_i - 2
\alpha \beta g_{i+1} + \alpha^2 g_{i+2}) + \dots + \\ + t^{r-i-2}
u^i &(\beta^2 g_{j-2} - 2 \alpha \beta g_{j-1} + \alpha^2 g_j) +
t^{r-i-1} u^{i-1} (\beta^2 g_{j-1} - \alpha \beta g_j)) = 0
\end{split} \end{equation*} and the claim follows.

Let $ (\alpha, \beta) = (0, 1).$ Then, as  explained before, the double structure $ X $ associated to $ \mu = (0, \dots, 0, 1) $ is contained in the cone over the rational normal curve $ C' $ of $ \pp^{r-1} \cong H = V(x_r) $ defined by the $ 2 \times 2 $ minors of the matrix $$ \left( \begin{array}{cccc} x_0 & x_1 & \dots & x_{r-2} \\ x_1 & x_2 & \dots & x_{r-1} \end{array} \right).$$ The hyperplane $ H $ intersects $ C $ in $ A (1: 0 : \dots : 0) $ with multiplicity $ r,$ and so $ X \cap H $ is supported on $ A $ and has degree $ 2r.$ Moreover, $ X \cap H $ is contained in the rational normal curve $ C' $ and so it is arithmetically Gorenstein. In fact, $ \deg(X \cap H) = 2 \deg(C') + 2,$ $ Pic(C') \cong \zz,$ and so $ X \cap H \in \vert 2H' - K' \vert,$ where $ H' $ is the class of a hyperplane section of $ C' $ and $ K' $ is the canonical divisor of $ C',$ and every divisor of the system $ dH' - K' $ is arithmetically Gorenstein, for every $ d \geq 0 $ (\cite{mig-book}, Theorem 4.2.8). Hence, $ X $ is arithmetically Cohen--Macaulay with an arithmetically Gorenstein hyperplane section, i.e. $ X $ is arithmetically Gorenstein, because the graded Betti numbers of the minimal free resolutions of $ I_X $ and $ I_{X \cap H \vert H} $ are the same. In fact, the only irreducible component of $ X $ is $ C $ that is non--degenerate, and so $ x_r $ is not a $ 0$--divisor for $ R/I_X.$ Hence, the hyperplane $ H = V(x_r) $ we considered is general enough for $ X $ to let the proof of (\cite{mig-book}, Theorem 1.3.6) work (see also \cite{mig-book}, Remark 1.3.9).

The family $ H(2r, r+1, r) $ is irreducible, and the arithmetically Gorenstein locus in it is not empty. From the semicontinuity of the Betti numbers in an irreducible family (\cite{greco-bor-2}) it follows that an open subscheme of $ H(2r, r+1, r) $ parametrizes arithmetically Gorenstein schemes and the claim follows.
\end{proof}

\begin{remark} \rm In the case we just studied, the degree and the
genus of $ X $ are the ones of a canonical curve in $ \pp^r,$ and
so the result we proved is not unexpected. In fact, in \cite{fong}, the author proved that $ H(2r, r+1, r) $ is contained in the closure of the component of the canonical curves.
\end{remark}

Now, we consider the second case, namely $ (2r,g,n) = (2r, 1,
2r-1).$

\begin{theorem} \label{existence} Let $ (2r,g,n) = (2r,1,2r-1).$ Then for a general
map $ \mu: \oo^{r-1}_{\pp^1}(-r-2) \oplus \oo^{r-1}_{\pp^1}(-r)
\to \oo_{\pp^1}(-2) $ we get a non--degenerate arithmetically
Gorenstein curve, where general means that $ \mu_{\vert}:
\oo^{r-1}_{\pp^1}(-r) \to \oo_{\pp^1}(-2) $ has no degree $ 0 $
syzygy.
\end{theorem}

\begin{proof} Assume that the restriction $ \mu_1 $ of $ \mu $ to
$ \oo_{\pp^1}^{r-1}(-r-2) $ is the null map.

If $ r = 2,$ then $ \mu = (0, 1) $ and hence $ X $ is defined by
the ideal $ I_X = (x_0x_2 - x_1^2, x_3^2).$ Then $ X $ is a
complete intersection of a cone and a double plane.

Assume now that $ r \geq 3,$ and furthermore assume that the
restriction $ \mu_2 $ of $ \mu $ to $ \oo^{r-1}_{\pp^1}(-r) $ is
given by $ \mu_2 = (t^{r-2}, t^{r-3} u, \dots, u^{r-2}).$ Of
course, $ \mu_2 $ has no degree $ 0 $ syzygy. The map $ \mu \circ
\psi $ is given by the matrix $ (0, \dots, 0, t^{r-2}, t^{r-3} u,
\dots, u^{r-2}).$ By using the procedure described in section
\ref{s-construction}, we get that the double structure $ X $ is
defined by the ideal $ I_X = I_{C,L} + I_L^2 + J $ where $ J =
\langle x_1 x_{r+1} - x_0 x_{r+2}, \dots, x_r x_{r+1} - x_{r-1}
x_{r+2}, \dots, x_1 x_{2r-2} - x_0 x_{2r-1}, \dots, x_r x_{2r-2} -
x_{r-1} x_{2r-1} \rangle.$ It is evident that the ideal $ I_S $
defined by the $ 2 \times 2 $ minors of the matrix $$ \left(
\begin{array}{cccccc} x_0 & \dots & x_{r-1} & x_{r+1} & \dots &
x_{2r-2} \\ x_1 & \dots & x_r & x_{r+2} & \dots & x_{2r-1}
\end{array} \right) $$ is contained in $ I_X,$ i.e. $ X $ is
contained in $ S $ which is a smooth rational normal scroll
surface $ \pp(\oo_{\pp^1} \oplus \oo_{\pp^1}(-2)) $ embedded via
the complete linear system $ \vert \xi + r f \vert $
(\cite{eisenbud-book}, exercise A2.22), where $ \xi $ is the class
of the rational normal curve of minimal degree $ r-2 $ contained
in $ S $ and $ f $ is a fibre. On $ S,$ we have that $ \xi^2 = -2, \xi \cdot f = 1, f^2 = 0.$ Moreover, the canonical divisor of $ S $ is $ K_S = -2 \xi -2 f,$ and the hyperplane section class is $ H = \xi + r f $ (\cite{hart}, Lemma 2.10). Then $ X \in \vert a \xi + b f
\vert,$ with $ a = 2, b = 4,$ by adjunction, and so $ X $ is an
anticanonical divisor on $ S $ and so it is a non--degenerate
arithmetically Gorenstein curve (\cite{mig-book}, Theorem 4.2.8). To complete the proof, we show
that the curve $ X $ we constructed before is the only double
structure on $ C $ of arithmetic genus $ 1,$ up to automorphisms of $ \pp^{2r-1},$
which is the content of next Theorem \ref{unique}.
\end{proof}

\begin{theorem} \label{unique} Let $ C \subset L \cong \pp^r \subset
\pp^{2r-1} $ be a rational normal curve of degree $ r.$ Then,
there exists only one non--degenerate double structure $ X $ on $
C $ of arithmetic genus $ 1,$ up to automorphisms of $ \pp^{2r-1}.$
\end{theorem}

\begin{proof} To make the proof more readable, we choose the
coordinates of $ \pp^{2r-1} $ as $ x_0, \dots, x_r,$ $ y_1, \dots,
y_{r-1},$ where $ L = V(y_1, \dots, y_{r-1}).$

As in the proof of the previous Theorem, let $ \mu_1 $ and $ \mu_2
$ be the restrictions of $ \mu $ to $ \oo^{r-1}_{\pp^1}(-r-2) $
and to $ \oo^{r-1}_{\pp^1}(-r),$ respectively. Assume first that $
\mu_2 = (l_1, \dots, l_{r-1}).$ The forms $ l_1, \dots, l_{r-1} $ are
linearly dependent if, and only if, they have a degree $ 0 $
syzygy, that, of course, is also a degree $ 0 $ syzygy of $ \mu
\circ \psi.$ So, $ X $ is degenerate if, and only if, $ l_1,
\dots, l_{r-1} $ are linearly dependent. Hence, we can assume that
$ l_1, \dots, l_{r-1} $ are linearly independent, and so there
exists an invertible matrix $ P \in GL_{r-1}(K) $ such that $$
\left( \begin{array}{c} l_1 \\ \vdots
\\ l_{r-1} \end{array} \right) = P \left( \begin{array}{c} t^{r-2}
\\ \vdots \\ u^{r-2} \end{array} \right).$$

Going back to the construction, it is clear that the choice of the
generators of $ L $ plays no role when we restrict the maps to $
\pp^1 $ and so, if we say that $ I_L $ is generated by $$ P^{-1}
\left( \begin{array}{c} y_1 \\ \vdots \\ y_{r-1} \end{array}
\right) $$ and we change bases in $ \oo^{r-1}_{\pp^1}(-r) $ by
using $ P^{-1} $ we get that $ I_L $ is generated by $ y_1, \dots,
y_{r-1} $ and $ \mu_2 = (t^{r-2}, \dots, u^{r-2}).$

Let $ \mu_1 = (p_1, \dots, p_{r-1}),$ where $ p_i = p_{i0} t^r +
p_{i1} t^{r-1}u + \dots + p_{ir} u^r.$ The map $ \mu \circ \psi $
has the first $ \binom r2 $ entries that are combinations of $
t^{r-2}, \dots, u^{r-2} $ with coefficients $ p_1, \dots, p_{r-1}
$ and the last $ r-1 $ entries which are equal to $ t^{r-2}, \dots
u^{r-2}.$ Hence, the syzygies of $ \mu \circ \psi $ can be easily
computed and we get that the defining ideal $ I_X $ of $ X $ is
generated by $ x_{i-1} x_{i+1} - x_i^2 - p_i y_i, i=1, \dots,
r-1,$ by $ x_i x_{j+1} - x_{i+1} x_j - y_{i+1} p_j - \dots - y_j
p_{i+1} , 0 \leq i < j-1 \leq r-2,$ by $ x_i y_j - x_{i-1}
y_{j+1}, i = 1, \dots, r, j = 1, \dots, r-2,$ and by $ y_i y_j, 1
\leq i \leq j \leq r-1,$ where, with abuse of notation, we set $
p_i $ also the only linear form in $ x_0, \dots, x_r $ that is
equal to $ p_i $ when restricted to $ \pp^1,$ i.e. $ p_i = p_{i0}
x_0 + \dots + p_{ir} x_r.$

We look for the required change of coordinates in the form $$ x_i
= z_i + a_{i1} y_1 + \dots a_{i,r-1} y_{r-1} \qquad i = 0, \dots,
r $$ and we fix the remaining variables $ y_1, \dots, y_{r-1}.$

Our goal is to prove that we can choose the $ a_{ij}$'s in such a
way that, in the new coordinate system, $ X $ is defined by the
ideal $ J $ generated by $ z_{i-1} z_{i+1} - z_i^2, z_i z_{j+1} -
z_{i+1} z_j, z_i y_j - z_{i-1} y_{j+1}, y_i y_j,$ where the indices
vary in the same ranges as before.

If we apply the change of coordinates to the last generators of $
I_X $ then they do not change, because the variables $ y_i, \dots
y_{r-1} $ are fixed. If we apply the change of coordinates to the
generators of the form $ x_i y_j - x_{i-1} y_{j+1},$ we get that $
z_i y_j - z_{i-1} y_{j+1} \in J $ because $ y_h y_k \in J,$ for $
i = 1, \dots, r, j = 1, \dots, r-2,$ and $ 1 \leq h \leq k \leq
r-1.$

By applying the change of coordinates to $ x_{i-1} x_{i+1} - x_i^2
- y_i p_i $ we get \begin{equation*} \begin{split} ( & z_{i-1} +
\sum_{j=1}^{r-1} a_{i-1,j} y_j) (z_{i+1} + \sum_{j=1}^{r-1}
a_{i+1,j} y_j) - (z_i + \sum_{j=1}^{r-1} a_{ij} y_j)^2 - y_i
\sum_{k=0}^r p_{ik}z_k = \\ = & z_{i-1} z_{i+1} - z_i^2 -
\sum_{j=1}^{r-1} a_{i-1,j} z_{i+1} y_j - \sum_{j=1}^{r-1}
a_{i+1,j} z_{i-1} y_j + 2 \sum_{j=1}^{r-1} a_{ij} z_i y_j -
\sum_{k=0}^r p_{ik} z_k y_i =^{(*)} \\
= & z_{i-1} z_{i+1} - z_i^2 - \sum_{j=1}^{r-1} a_{i-1,j}
z_{i+j+2-r} y_{r-1} - \sum_{j=1}^{r-1} a_{i+1,j} z_{i+j-r} y_{r-1}
+ \\ & + 2 \sum_{j=1}^{r-1} a_{ij} z_{i+j+1-r} y_{r-1} -
\sum_{k=0}^r p_{ik} z_{k+i+1-r} y_{r-1}, \end{split}
\end{equation*}
where, in $ (*),$ we use the fact that $ z_i y_j - z_{i-1} y_{j+1}
\in J $ and the convention that we can use $ z_x y_{r-1} $ with $
x < 0 $ to mean $ z_0 y_{r+x-1}.$

Hence, we get the following linear equations in the $ a_{ij}$'s:
\begin{equation} \label{main-rel} a_{i-1,h-1} - 2 a_{ih} +
a_{i+1,h+1} - p_{ih} = 0 \end{equation} for $ h = 0, \dots, r $
where we assume that $ a_{ij} = 0 $ if $ j \leq 0 $ or $ j \geq
r,$ for whatever $ i.$

With analogous computations, if we apply the change of coordinates
to $ x_i x_{j+1} - x_{i+1} x_j - y_{i+1} p_j - \dots - y_j p_{i+1}
$ we get the following linear equations $$ a_{i,m-j-2+r} +
a_{j+1,m-i-1+r} - a_{j,m-i-2+r} - a_{i+1,m-j-1+r} - \sum_{t=i+1}^j
p_{i+j+1-t,m-t-1+r} = 0.$$

It is an easy computation to show that those last equations depend
linearly from the previous ones. For example, if we subtract from
the last equation the one among (\ref{main-rel}) we get setting $
i = j, h = m-i-2+r,$ we have the relation $$ a_{i,m-j-2+r} +
a_{j,m-i-2+r} - a_{j-1,m-i-3+r} - a_{i+1,m-j-1+r} - \sum_{t=i+2}^j
p_{i+j+1-t,m-t-1+r} = 0,$$ which is again of the same form, but
with smaller difference between the first subscripts. By
iterating, we get that all of them linearly depend from the
equations (\ref{main-rel}). \medskip

Now, we prove that the linear system (\ref{main-rel}) has one
solution. To this aim, we collect the equation according to the
difference $ i-h $ of the subscripts of the variables involved. In
fact, notice that in each equation, the difference is constant.

At first, assume that the difference is $ i-h = 0.$ Then, we get
the following linear system $$ \left\{ \begin{array}{l} -2 a_{11}
+ a_{22} = p_{11} \\ a_{11} - 2 a_{22} + a_{33} = p_{22} \\ \vdots \\
a_{i-1,1-i} - 2 a_{ii} + a_{i+1,i+1} = p_{ii} \\ \vdots \\
a_{r-2,r-2} - 2 a_{r-1,r-1} = p_{r-1,r-1} \end{array} \right.$$
The coefficient matrix $ M_{r-1} = (m_{ij}) $ has entries equal to
$$ m_{ij} = \left\{ \begin{array}{cl} -2 & \mbox{if } i = j \\ 1 &
\mbox{if } \vert i-j \vert = 1 \\ 0 & \mbox{otherwise} \end{array}
\right. $$ and its determinant is equal to $ 2p+1 $ if $ r-1=2p,$
or to $ -2p $ if $ r-1 = 2p-1.$ In fact, by the Laplace formula,
$$ \det(M_{r-1}) = -2 \det(M_{r-2}) - \det(M_{r-3}),$$ by direct
computation $ \det(M_1) = -2, \det(M_2) = 3 $ and the claim can be
easily proved by induction. Hence, the previous linear system has
one solution, by Cramer's rule.

Assume now that the difference is equal to $ i-h = k > 0.$ Hence,
the corresponding linear system is $$ \left\{ \begin{array}{l}
a_{k+1,1}  = p_{k0} \\ -2 a_{k+1,1} + a_{k+2,2} = p_{k+1,1} \\
a_{k+1,1} - 2 a_{k+2,2} + a_{k+3,3} = p_{k+2,2} \\ \vdots \\
a_{r-2, r-2-k} - 2 a_{r-1,r-1-k} + a_{r,r-k} = p_{r-2,r-2-k}
\end{array} \right. $$ and it has one solution for every $ k.$

Analogously, the system with $ i-h = k < 0 $ has one solution for
every $ k $ and the claim follows.
\end{proof}

\begin{remark} \rm Let $ X $ be the double structure on a rational normal curve defined in the proof of Theorem \ref{existence}. There is a natural map $ \Psi: Aut(\pp^{2r-1}) \to H(2r, 1, 2r-1) $ defined as $ \Psi(g) = g(X) $ where $ g(X) $ is the double structure we get by applying $ g $ to $ X.$ Previous Theorem \ref{unique} is equivalent to $ \ker \Psi = Aut(C).$ In fact, every automorphism $ g $ of $ C $ extends to an automorphism $ g' $ of $ L $ that fixes $ C.$ $ g' $ can be further extended to an automorphism $ g'' $ of $ \pp^{2r-1} $ that fixes $ L.$ For such a $ g'' $ we have that $ g''(X) = X.$ Hence, $ Aut(C) \subseteq \ker \Psi.$ By a dimension count, we get that $ \Psi $ is surjective if, and only if, $ \ker \Psi = Aut(C).$
\end{remark}

Now, we apply the previous results to Gorenstein liaison.
\begin{corollary} A rational normal curve $ C \subseteq L \cong \pp^r \subseteq \pp^n $ of degree $ r $ is self-linked if, and only if, either $ n = r $ or $ n = 2r-1.$
\end{corollary}

\begin{proof} $ C $ is self-linked if, and only if, there exists a double structure supported on $ C $ that is arithmetically Gorenstein. The claim is then a direct consequence of Theorems \ref{canonical} and \ref{existence}.
\end{proof}


\section{Double conics}

In this section, we prove that the general double structure of genus $ g
\leq -2 $ supported on a smooth conic is a smooth point in its Hilbert
scheme. Moreover, they are the general element of an
irreducible component in the same range of the arithmetic genus. On the other hand, if such double structures are
contained in $ \pp^3,$ and their genus satisfies $ g \geq -1,$
then we identify the general element of the irreducible component
containing the considered double structures. To achieve the
result, we compare the dimension of the family of the double
structures of fixed genus with the dimension of $ H^0(X, \ns_X),$
global sections of the normal sheaf of a suitable double conic $
X.$ In fact, it is well known that $ H^0(X, \ns_X) $ can be
identified with the tangent space to the Hilbert scheme at $ X.$
To get the desired results, we consider first a suitable double
conic $ X \subset \pp^3.$ We describe its ideal $ I_X $ and the
minimal free resolution $$ 0 \to F_3 \to F_2 \to F_1 \to I_X \to
0,$$ where $ F_i $ is a free $ R = K[x,y,z,w]-$module, then the $
R-$module structure of the global sections of $ H^0_*(X, \oo_X) =
\bigoplus_{j \in \zz} H^0(X, \oo_X(j)) $ and finally we compute $
H^0(X, \ns_X) $ as the degree $ 0 $ elements of $$ \ker(\Hom(F_1,
H^0_*(X, \oo_X)) \to \Hom(F_2, H^0_*(X, \oo_X))).$$ The result for
a general double conic $ X \subset \pp^n $ follows from the
smoothness of the Hilbert scheme at a degenerate double conic.
Furthermore, in $ \pp^3,$ we distinguish the case $ g(X) $ odd
from the case $ g(X) $ even, because the ideals have a different
minimal number of generators, and so their minimal free
resolutions have not comparable free modules and maps. Of course,
even if there are differences, we use the same arguments in both
cases.

\subsection{Case $ g $ odd, i.e. $ a = 2b $} In this subsection, we will
use the following running notation. We set $ R = K[x,y,z,w],$ and $ C \subset
\pp^3 = \Proj(R) $ is the conic defined by the ideal $ I_C =
(xz-y^2, w).$ Let $ j: \pp^1 \to C $ be the isomorphism defined as
$ j(t:u) = (t^2: tu: u^2: 0).$ Finally, we set $ \mu: \oo_{\pp^1}(-4) \oplus
\oo_{\pp^1}(-2) \to \oo_{\pp^1}(-4+2b) $ to be the map defined as $
\mu = (u^{2b}, t^{2b-2}).$

\begin{proposition} \label{ideal} If $ X $ is the doubling of $ C $
associated to $ \mu,$ then \begin{enumerate} \item $ I_X = \langle
w^2, w(xz-y^2), (xz-y^2)^2, x^{b-1} (xz-y^2) - z^b w \rangle;$
\item $ w^2, w(xz-y^2), (xz-y^2)^2, x^{b-1} (xz-y^2) - z^b w $ is
a Gr\"obner basis of $ I_X $ with respect to the reverse
lexicographic order; \item the minimal free resolution of $ I_X
$ is
$$ 0 \to R(-b-4)
\stackrel{\delta_2}{\lra} \begin{array}{c} R(-4) \\ \oplus \\
R(-5) \\ \oplus \\ R(-b-2) \\ \oplus \\ R(-b-3) \end{array}
\stackrel{\delta_1}{\lra} \begin{array}{c} R(-2) \\ \oplus \\
R(-3) \\ \oplus \\ R(-4) \\ \oplus \\ R(-b-1) \end{array} \lra I_X
\to 0,$$ where the maps $ \delta_1 $ and $ \delta_2 $ will be
described in the proof.
\end{enumerate}
\end{proposition}

\begin{proof} The syzygies of $ H^0_*(j_* \mu) $ are generated by $ M = \left( \begin{array}{c} x^{b-1} \\ -z^b \end{array} \right) $ and so the saturated ideal $ I_X $ of $ X $
is generated by $ I_C^2 + [I_C] M,$ that is to say, $$ I_X =
\langle w^2, w(xz-y^2), (xz-y^2)^2, x^{b-1} (xz-y^2) - z^b w
\rangle.$$  The generators are a Gr\"obner basis of $ I_X $
because their $ S-$polynomials reduces to $ 0 $ via themselves
(\cite{eisenbud-book}, Theorem 15.8). Moreover, the free $
R-$module $ F_1 $ follows. Let $ (g_1, \dots, g_4) $ be a syzygy
of $ I_X.$ Then, in $ R,$ we have
$$ w^2 g_1 + w(xz-y^2) g_2 + (xz-y^2)^2 g_3 + x^{b-1} (xz-y^2) g_4
- z^b w g_4 = 0,$$ that can be rewritten as $ w ( w g_1 + (xz-y^2)
g_2 - z^b g_4) + (xz-y^2) ((xz-y^2) g_3 + x^{b-1} g_4) = 0.$ The
two polynomials $ w, xz-y^2 $ are a regular sequence, and so there
exists $ g \in R $ such that $$ \begin{array}{l}  w g + (xz-y^2)
g_3 + x^{b-1} g_4 = 0 \\ w g_1 + (xz-y^2) (g_2 -g) - z^b g_4 = 0.
\end{array} $$ Both $ w, xz-y^2, x^{b-1} $ and $ w, xz-y^2, z^b $
form a regular sequence, and so we have $$ \left( \begin{array}{c} g
\\ g_3 \\ g_4 \end{array} \right) = \left( \begin{array}{ccc} 0 &
x^{b-1} & -(xz-y^2) \\ -x^{b-1} & 0 & w \\ xz-y^2 & -w & 0
\end{array} \right) \left( \begin{array}{c} f_1 \\ f_2 \\ f_3
\end{array} \right) $$ and $$ \left( \begin{array}{c} g_1 \\ g_2-g
\\ -g_4 \end{array} \right) = \left( \begin{array}{ccc} 0 & z^b &
-(xz-y^2) \\ -z^b & 0 & w \\ xz-y^2 & -w & 0 \end{array} \right)
\left( \begin{array}{c} f_4 \\ f_5 \\ f_6 \end{array} \right).$$
By comparing the value of $ g_4 $ from the two expressions above, we get
the equation $ (xz-y^2) (f_1 + f_4) = w (f_2 + f_5).$ By using the
same argument as before, there exists $ h \in R $ such that $ f_4
= -f_1 + w h, f_5 = -f_2 + (xz-y^2) h.$ Hence, it holds $$ \left(
\begin{array}{c} g_1 \\ g_2 \\ g_3 \\ g_4 \end{array} \right) =
\left( \begin{array}{cccc} -(xz-y^2) & 0 & -z^b & 0 \\ w &
-(xz-y^2) & x^{b-1} & z^b \\ 0 & w & 0 & -x^{b-1} \\ 0 & 0 & -w &
xz-y^2 \end{array} \right) \left( \begin{array}{c} f_6 - z^b h \\
f_3 \\ f_2 \\ f_1 \end{array} \right),$$ and the $ 4 \times 4 $
matrix represents the map $ \delta_1.$ Of course, the free $
R-$module $ F_2 $ follows from $ F_1 $ and from the degrees of the
entries of the map $ \delta_1.$

The second syzygies of $ I_X $ can be computed as the first ones,
and we get $$ \delta_2 = \left( \begin{array}{c} -z^b \\ x^{b-1}
\\ xz-y^2 \\ -w \end{array} \right).$$
\end{proof}

The last statement of next Proposition is due to the anonymous
referee that we thank once more.

\begin{proposition} $ X $ has
genus $ g(X) = 3-2b,$ and the Hartshorne-Rao function of $ X $ is
$$ h^1 \ii_X(j) = \left\{ \begin{array}{cl} 2(j+b)-3 & \mbox{if }
-b+2 \leq j \leq 0 \\ 2b-2 & \mbox{if } j = 1 \\ 2(b-j)+1 &
\mbox{if } 2 \leq j \leq b \\ 0 & \mbox{otherwise}. \end{array}
\right. $$ Moreover, the Hartshorne-Rao module of $ X $ is
isomorphic to $$ R/\langle w, xz-y^2, x^{b-1}, z^{b}
\rangle(-b+2).$$
\end{proposition}

\begin{proof} The genus and the Hartshorne--Rao function can be
computed by using results from section \ref{properties-of-X}.

For the last statement, we first remark that double conics are
minimal curves in the sense of \cite{mdp}. Otherwise, a double conic would be bilinked down to a degree two curve.
But the Hartshorne--Rao function of a degree two curve increases by
at most one, and so we can exclude the case. For minimal curves,
the map $ \delta_2^\vee $ begins a minimal free resolution of $
H^1_* \ii_X $ and so the claim follows because the entries of $
\delta_2 $ are a regular sequence.
\end{proof}

Now, we can compute the dimension of the degree $ d $ global
sections of the structure sheaf of $ X.$ In fact, it holds
\begin{proposition} \label{h0OX} $ h^0(X, \oo_X(d)) = 4d+2b-2 $ if $ d \geq 2.$
\end{proposition}

\begin{proof} The short exact sequence $$ 0 \to \frac{R}{I_X}
\to H^0_*(X, \oo_X) \to H^1_* \ii_X \to 0 $$ allows us to prove
the result.
\end{proof}

Now, we describe the elements of $ H^0(X, \oo_X(d)) $ for every $
d \geq 2.$ To start, we can easily describe the elements of $
H^0(X, \oo_X(d)) $ for $ d \geq b+1,$ because, in the considered
range, we have $ (R/I_X)_d = H^0(X, \oo_X(d)),$ and so it holds
\begin{proposition} Let $ d \geq b+1.$ Then, $ H^0(X, \oo_X(d)) = V_d $
where $$ V_d = \{ p_1 + y p_2 + (xz-y^2) p_3 + (xz-y^2) y p_4 + w
x^{d-b} p_5 + w x^{d-b-1} y p_6 \vert p_i \in K[x,z] \} $$ and the
degrees of the $ p_i'$s are fixed in such a way that the elements
in $ V_d $ are homogeneous of degree $ d.$
\end{proposition}

\begin{proof} The generators of $ I_X $ are a Gr\"obner basis and
so the initial ideal $ in(I_X) $ of $ I_X $ with respect to the
reverse lexicographic order is generated by $ w^2, y^2 w, y^4,
x^{b-1} y^2.$ Hence, the elements in $ V_d $ are in normal form
with respect to $ I_X $ and so they are linearly independent.
\end{proof}

To describe the elements of $ H^0(X, \oo_X(d)) $ for $ 2 \leq d
\leq b,$ we first define a suitable global section $ \xi $ of
degree $ -b+2,$ and then we compute all the global sections by
using $ \xi $ and the elements in $ (R/I_X)_d.$

\begin{definition} Let $ \xi \in H^0(X, \oo_X(-b+2)) $ be the
global section of $ X $ defined as $$ \xi = \frac{w}{x^{b-1}} =
\frac{xz-y^2}{z^b}.$$
\end{definition}

The global section $ \xi $ is well defined because for no closed point on $ C $ both
$ x $ and $ z $ can be equal to $ 0,$ and because the two
descriptions agree on the overlap (in fact $ x^{b-1} (xz-y^2) -
z^b w \in I_X$).

\begin{proposition} $ \xi $ verifies the following equalities
\begin{enumerate} \item $ w \xi = (xz-y^2) \xi = 0;$ \item $ x^{b-1}
\xi = w;$ \item $ z^b \xi = xz-y^2.$ \end{enumerate}
\end{proposition}

\begin{proof} They follow easily from the definition of $ \xi $ and
from the knowledge of the ideal $ I_X.$
\end{proof}

\begin{proposition} \label{g-s-1} Let $ 2 \leq d
\leq b.$ Then, $ H^0(X, \oo_X(d)) = V_d $ where \begin{equation*}
\begin{split} V_d = & \{p_1 + y p_2 + (xz-y^2) p_3 + (xz-y^2) y p_4 + w p_5
+ wy p_6 + \\ & \xi (x^{d-1} z^d q_1 + x^{d-2} y z^{d-1} q_2)
\vert p_i \in K[x,z], q_j \in K[x,z] \} \end{split}
\end{equation*} and the elements in $ V_d $ are homogeneous
of degree $ d.$
\end{proposition}

\begin{proof} Of course, $ V_d \subseteq H^0(X, \oo_X(d)).$

As before, the elements $ p_1 + y p_2 + (xz-y^2) p_3 + (xz-y^2) y
p_4 + w p_5 + wy p_6 $ are in normal form with respect to $ I_X $
and so they are linearly independent. Let $ \pi: H^0_* \oo_X \to
H^1_* \ii_X.$ It is evident that $ \xi \notin \ker(\pi) = R/I_X $
and so $ \pi(\xi) \not= 0.$ Hence $ \pi( \xi(x^{d-1} z^d q_1 +
x^{d-2} y z^{d-1} q_2)) = \pi(\xi) (x^{d-1} z^d q_1 + x^{d-2} y
z^{d-1} q_2) \in R/\langle w, xz-y^2, x^{b-1}, z^b\rangle (-b+2).$
The generators of $ \langle w, xz-y^2, x^{b-1}, z^b\rangle $ are a
Gr\"obner basis and $ x^{d-1} z^d q_1 + x^{d-2} y z^{d-1} q_2 $
are in normal form with respect to the given Gr\"obner basis.
Hence, they are linearly independent, and so $ V_d $ has dimension
$ \dim V_d = 4d + 2b -2,$ and the equality $ H^0(X, \oo_X(d)) =
V_d $ holds.
\end{proof}

Now, we compute the degree $ 0 $ global sections of $ H^0(X,
\ns_X) $ as $$ H^0(X, \ns_X) = \ker(\Hom(F_1, H^0_*(X, \oo_X))
\stackrel{\delta_1^\vee}{\lra} \Hom(F_2, H^0(X, \oo_X)))_0,$$
where, if $ \varphi \in \Hom(F_1, H^0_*(X, \oo_X)) $ then $
\delta_1^\vee(\varphi) = \varphi \circ \delta_1.$

Let $ F_1 = \oplus_{i=1}^4 R e_i $ with $ \deg(e_1) = 2, \deg(e_2)
= 3, \deg(e_3) = 4, \deg(e_4) = b+1,$ and assume $ b \geq 4.$ $
\varphi \in \ker(\delta_1^\vee) $ if, and only if, the following
system is satisfied: \begin{equation} \label{system} \left\{
\begin{array}{l} (xz-y^2) \varphi(e_1) - w
\varphi(e_2) = 0 \\ (xz-y^2) \varphi(e_2) - w \varphi(e_3) = 0 \\
z^b \varphi(e_1) - x^{b-1} \varphi(e_2) + w \varphi(e_4) = 0 \\
z^b \varphi(e_2) - x^{b-1} \varphi(e_3) + (xz-y^2) \varphi(e_4) =
0. \end{array} \right. \end{equation} with $ \varphi(e_i) \in
H^0(X, \oo_X(\deg(e_i))).$

To solve the system, we set \begin{itemize} \item[$\bullet$] $
\varphi(e_1) = p_{11} + y p_{12} + (xz-y^2) p_{13} + w p_{15} + wy
p_{16} + \xi (xz^2 q_{11} + yz q_{12}) $ with $ p_{1i} \in K[x,z],
q_{1j} \in K[x,z] $ and $ \deg(p_{11}) = 2, \deg(p_{12}) =
\deg(p_{15}) = 1, \deg(p_{13}) = \deg(p_{16}) = 0, \deg(q_{11}) =
b-3, \deg(q_{12}) = b-2;$ \item[$\bullet$] $ \varphi(e_2) = p_{21}
+ y p_{22} + (xz-y^2) p_{23} + (xz-y^2) y p_{24} + w p_{25} + wy
p_{26} + \xi (x^2z^3 q_{21} + xyz^2 q_{22}) $ with $ p_{2i} \in
K[x,z], q_{2j} \in K[x,z] $ and $ \deg(p_{21}) = 3, \deg(p_{22}) =
\deg(p_{25}) = 2, \deg(p_{23}) = \deg(p_{26}) = 1, \deg(p_{24}) =
0, \deg(q_{21}) = b-4, \deg(q_{22}) = b-3;$ \item[$\bullet$] $
\varphi(e_3) = p_{31} + y p_{32} + (xz-y^2) p_{33} + (xz-y^2) y
p_{34} + w p_{35} + wy p_{36} + \xi (x^3z^4 q_{31} + x^2yz^3
q_{32}) $ with $ p_{3i} \in K[x,z], q_{3j} \in K[x,z] $ and $
\deg(p_{31}) = 4, \deg(p_{32}) = \deg(p_{35}) = 3, \deg(p_{33}) =
\deg(p_{36}) = 2, \deg(p_{34}) = 1, \deg(q_{31}) = b-5,
\deg(q_{32}) = b-4;$ \item[$\bullet$] $ \varphi(e_4) = p_{41} + y
p_{42} + (xz-y^2) p_{43} + (xz-y^2) y p_{44} + w x p_{45} + wy
p_{46} $ with $ p_{4i} \in K[x,z] $ and $ \deg(p_{41}) = b+1,
\deg(p_{42}) = b, \deg(p_{45}) = \deg(p_{43}) = \deg(p_{46}) =
b-1, \deg(p_{44}) = b-2.$ \end{itemize}

\begin{theorem} \label{thm1} With the notation as above, $ h^0(X, \ns_X) = 7
+ 4b = 13 - 2g,$ where $ g $ is the arithmetic genus of $ X.$
\end{theorem}

\begin{proof} To get the claim, we have to compute the elements in $ H^0(X, \ns_X).$

\noindent \em Claim: \rm
$ \varphi \in H^0(X,
\ns_X) $ if, and only if, $ \varphi(e_1) = 2w P_1 + 2wy P_2,
\varphi(e_2) = (xz-y^2) P_1 + (xz-y^2) y P_2 + w P_3 + wy P_4,
\varphi(e_3) = 2(xz-y^2) P_3 + 2(xz-y^2) y P_4, \varphi(e_4) =
(x^{b-1} P_3 - z^b P_1) + y (x^{b-1} P_4 - z^b P_2) + (xz-y^2) P_5
+ (xz-y^2) y P_6 + w x P_7 + w y P_8 $ with $ P_i \in K[x,z] $ for
every $ i $ and $ \deg(P_1) = \deg(P_4) = 1, \deg(P_2) = 0,
\deg(P_3) = 2, \deg(P_5) = \deg(P_7) = \deg(P_8) = b-1, \deg(P_6)
= b-2.$

If the claim holds, then we get the dimension of $ H^0(X, \ns_X) $ with an easy parameter count. Hence, we prove the Claim.

It is easy to check that if $ \varphi $ satisfies
the given conditions, then $ \varphi \in H^0(X, \ns_X).$ Conversely, we solve the equations of the system
(\ref{system}) one at a time.

The first equation becomes $$ (xz-y^2) (p_{11} + y p_{12}) - w
(p_{21} + y p_{22}) = 0.$$ As a $ R/I_X-$module, the first syzygy
module of $ (xz-y^2, -w) $ is generated by $$ \left(
\begin{array}{ccccc} xz-y^2 & w & 0 & 0 & x^{b-1} \\ 0 & 0 & w &
xz-y^2 & z^b \end{array} \right) $$ and so we get the two
equations $$ p_{11} + y p_{12} = 0 \qquad p_{21} + y p_{22} = 0 $$
because $ p_{ij} \in K[x,z] $ and by degree argument.

Again as $ R/I_X-$module, the first syzygy module of $ (1, y) $ is
generated by $$ \left( \begin{array}{c} y  \\ -1 \end{array} \right)
$$ and so the solutions of two equations are $ p_{11} = p_{12} =
p_{21} = p_{22} = 0.$

The second equation, after substituting the computed solutions of
the first one, becomes $$ w(p_{31} + y p_{32}) = 0.$$ We have that
$ p_{3i} \in K[x,z] $ and so $ p_{31} + y p_{32} = 0.$ Because of
the same argument, and the knowledge of the first syzygy module of
$ (1, y),$ we get that the solutions of this equation are $ p_{31}
= p_{32} = 0.$

The third equation of the system (\ref{system}) becomes
\begin{equation*} \begin{split} (xz-y^2) & (x^{b-1} (p_{15} -
p_{23}) + x^{b-1} y (p_{16}-p_{24}) + z^b p_{13} + xz^2 q_{11} +
yz q_{12}) \\ & - w (x^{b-1} p_{25} + x^{b-1} y p_{26} - p_{41} -
y p_{42} + x^2z^3 q_{21} + xyz^2 q_{22}) = 0. \end{split}
\end{equation*} Because of the knowledge of the first syzygy
module of $ (xz-y^2, -w),$ from the previous equation we get the
two ones \begin{equation} x^{b-1} (p_{15} - p_{23}) + x^{b-1} y
(p_{16}-p_{24}) + z^b p_{13} + xz^2 q_{11} + yz q_{12} = x^{b-1}
(r_1 + y r_2) \end{equation} and \begin{equation} x^{b-1} p_{25} +
x^{b-1} y p_{26} - p_{41} - y p_{42} + x^2z^3 q_{21} + xyz^2
q_{22} = z^b (r_1 + y r_2) \end{equation} where $ r_1 \in
K[x,z]_1, r_2 \in K[x,z]_0.$

The first one can be rewritten as $$ x^{b-1} (p_{15} - p_{23} -
r_1) + z^b p_{13} + xz^2 q_{11} + y (x^{b-1}(p_{16}-p_{24} - r_2)
+ z q_{12}) = 0.$$ From the knowledge of the first syzygy module
of $ (1, y) $ we get $$ x^{b-1} (p_{15} - p_{23} - r_1)
+ z^b p_{13} + xz^2 q_{11} = 0 $$ and $$ x^{b-1}(p_{16}-p_{24} -
r_2) + z q_{12} = 0.$$ Hence, $ p_{15} = p_{23} + r_1, p_{13} = 0,
q_{11} = 0, p_{16} = p_{24} + r_2, q_{12} = 0.$

The second equation can be rewritten as $$ x^{b-1} p_{25} - p_{41}
+ x^2z^3 q_{21} - z^b r_1 + y (- p_{42} + x^{b-1} p_{26} + xz^2
q_{22} - z^b r_2) = 0.$$ By using the same argument as before, we
get $$ x^{b-1} p_{25} - p_{41} + x^2z^3 q_{21} - z^b r_1
= 0 $$ and $$ - p_{42} + x^{b-1} p_{26} + xz^2 q_{22} - z^b r_2 =
0.$$ Hence, $ p_{41} = x^{b-1} p_{25} + x^2z^3 q_{21} - z^b r_1,
p_{42} = x^{b-1} p_{26} + xz^2 q_{22} - z^b r_2.$

The last equation of the system (\ref{system}) becomes
\begin{equation*} \begin{split} & (xz-y^2) (z^b(p_{23} - r_1) +
2x^2z^3 q_{21} + y (z^b(p_{24} - r_2) + 2xz^2 q_{22})) - \\ - w
&(z^b(p_{33} - 2 p_{25}) + x^{b-1} p_{35} + x^3z^4 q_{31} + y (
z^b (p_{34} - 2 p_{26}) + x^{b-1} p_{36} + x^2z^3 q_{32})) = 0.
\end{split} \end{equation*} Then, there exist $ r_3 \in K[x,z]_2,
r_4 \in K[x,z]_1 $ such that the two following equalities hold
\begin{equation} z^b (p_{23} - r_1) + 2 x^2 z^3 q_{21} + y
(z^b(p_{24} - r_2) + 2xz^2 q_{22}) = x^{b-1}(r_3 + y r_4)
\end{equation} and \begin{equation} z^b(p_{33} - 2 p_{25}) +
x^{b-1} p_{35} + x^3z^4 q_{31} + y ( z^b (p_{34} - 2 p_{26}) +
x^{b-1} p_{36} + x^2z^3 q_{32}) = z^b (r_3 + y r_4).
\end{equation} From the first equation, we get $$ z^b
(p_{23} - r_1) - x^{b-1} r_3 + 2x^2z^3 q_{21} = 0 $$ and $$ z^b
(p_{24} - r_2) - x^{b-1} r_4 + 2 x z^2 q_{22} = 0.$$ Hence, $ r_1
= p_{23}, r_2 = p_{24}, r_3 = r_4 = 0, q_{21} = q_{22} = 0.$ From
the second equation, we get $$  z^b(p_{33} - 2 p_{25}) + x^{b-1}
p_{35} + x^3z^4 q_{31} = 0 $$ and $$ z^b (p_{34} - 2 p_{26}) +
x^{b-1} p_{36} + x^2z^3 q_{32} = 0.$$ As before, we deduce that $
p_{33} = 2 p_{25}, p_{34} = 2 p_{26}, p_{35} = p_{36} = 0, q_{31}
= q_{32} = 0.$

Summarizing, we obtain the following: \begin{enumerate} \item $
\varphi(e_1) = 2w (p_{23} + y p_{24}),$ \item $ \varphi(e_2) =
(xz-y^2) (p_{23} + y p_{24}) + w (p_{25} + y p_{26});$ \item $
\varphi(e_3) = 2(xz-y^2) (p_{25} + y p_{26}),$ \item $\varphi(e_4)
= x^{b-1} (p_{25} + y p_{26}) - z^b (p_{23} + y p_{24}) + (xz-y^2)
(p_{43} + y p_{44}) + w (x p_{45} + y p_{46}),$ \end{enumerate}
and the claim follows with the obvious substitutions.
\end{proof}

\begin{remark} \rm We computed $ h^0(X, \ns_X) $ by using the function {\tt <normal\_sheaf} of the
computer algebra software Macaulay (see \cite{mac}), in the cases $ 0 \leq b \leq 3.$

If $ b = 3,$ or equivalently $ g = -3,$ we get that the dimension
of the degree $ 0 $ global sections of the normal sheaf of $ X $
with saturated ideal
$$ I_X = \langle w^2, w(xz-y^2), (xz-y^2)^2, x^2(xz-y^2)-z^3w
\rangle $$ is equal to $ h^0(X, \ns_X) = 19 = 13-2g.$

If $ b = 2,$ i.e. $ g = -1,$ the dimension of the degree $ 0 $
global sections of the normal sheaf of the double conic $ X $
defined by the ideal $$ I_X = \langle w^2, w(xz-y^2), (xz-y^2)^2,
x(xz-y^2)-z^2w \rangle $$ is equal to $ h^0(X, \ns_X) = 16 \not=
13-2g.$

If $ b = 1,$ i.e. $ g = 1,$ the double conic $ X $ is defined by
the ideal $$ I_X = \langle w^2, w(xz-y^2), (xz-y^2)^2, (xz-y^2)-z
w \rangle $$ and $ h^0(X, \ns_X) = 16 \not= 13-2g.$

At last, if $ b = 0,$ i.e. $ g = 3,$ then the double conic $ X $
is defined by the ideal $$ I_X = \langle w, (xz-y^2)^2 \rangle $$
and $ h^0(X, \ns_X) = 17 \not= 13-2g.$
\end{remark}

\subsection{Case $ g $ even, i.e. $ a = 2b+1 $}

In this subsection, we repeat what we did in the previous
subsection, by sketching the main differences.

The running notation of the subsection are the following. As
before, we set $ R = K[x,y,z,w] $ and $ C \subset \pp^3 = \Proj(R)
$ is the conic defined by the ideal $ I_C = (xz-y^2, w).$ $ C $ is
isomorphic to $ \pp^1 $ via $ j: \pp^1 \to C $ defined as $ j(t:u)
= (t^2: tu: u^2: 0).$ We set $ \mu: \oo_{\pp^1}(-4) \oplus
\oo_{\pp^1}(-2) \to \oo_{\pp^1}(-3+2b) $ defined as $ \mu =
(u^{2b+1}, t^{2b-1}).$

\begin{proposition} If $ X $ is the doubling of $ C $
associated to $ \mu,$ then \begin{enumerate} \item $ I_X = \langle
w^2, w(xz-y^2), (xz-y^2)^2, x^b (xz-y^2) - y z^b w, x^{b-1} y
(xz-y^2) - z^{b+1} w \rangle;$ \item the generators of $ I_X $ are
a Gr\"obner basis with respect to the reverse lexicographic order;
\item the minimal free resolution of $ I_X $ is
$$ 0 \to R^2(-b-4)
\stackrel{\delta_2}{\lra} \begin{array}{c} R(-4) \\ \oplus \\
R(-5) \\ \oplus \\ R^4(-b-3) \end{array}
\stackrel{\delta_1}{\lra} \begin{array}{c} R(-2) \\ \oplus \\
R(-3) \\ \oplus \\ R(-4) \\ \oplus \\ R^2(-b-2) \end{array} \lra
I_X \to 0,$$ where the maps $ \delta_1 $ and $ \delta_2 $ will be
described in the proof.
\end{enumerate}
\end{proposition}

\begin{proof} The first syzygy module of $ H^0_*(j_* \mu) = (yz^b, x^b) $ is generated by
$$ N = \left( \begin{array}{cc} x^b & x^{b-1} y \\ - y
z^b & -z^{b+1} \end{array} \right).$$ Let $ M $ be the matrix we get by reading $ N $ over $ R $ and not over $ R/I_C.$ Then, the saturated ideal of $ X
$ is $$ I_X = \langle w^2, w(xz-y^2), (xz-y^2)^2, x^b (xz-y^2) - y
z^b w, x^{b-1} y (xz-y^2) - z^{b+1} w \rangle.$$ The check on $
S-$polynomials holds on the generators of $ I_X $ and so they are
a Gr\"obner basis.

To compute the first syzygy module of $ I_X $ we proceed as in the
proof of Proposition \ref{ideal}. The computation is quite similar
and uses the same ideas. Hence, we write only the maps $ \delta_1
$ and $ \delta_2:$ $$ \delta_1 = \left(
\begin{array}{cccccc} xz-y^2 & 0 & -yz^b & 0 & 0 & z^{b+1} \\ -w
& xz-y^2 & x^b & z^b & 0 & -x^{b-1} y \\ 0 & -w & 0 & 0 & -x^{b-1}
& 0 \\ 0 & 0 & -w & -y & z & 0 \\ 0 & 0 & 0 & x & -y & w
\end{array} \right) $$ and $$ \delta_2 = \left( \begin{array}{cc}
z^b & 0 \\ 0 & x^{b-1} \\ -y & -z \\ w & 0 \\ 0 & -w \\ -x & -y
\end{array} \right).$$
\end{proof}

\begin{proposition} $ X $ has
genus $ g(X) = 2-2b,$ the Hartshorne-Rao function of $ X $ is
equal to $$ h^1 \ii_X(j) = \left\{ \begin{array}{cl} 2(b+j-1) &
\mbox{ if } -b+1 \leq j \leq 0 \\ 2b-1 & \mbox{ if } j = 1 \\
2(b-j+1) & \mbox{ if } 2 \leq j \leq b+1 \\ 0 & \mbox{ otherwise}
\end{array} \right.$$ and $ \delta_2^\vee $ is a presentation
matrix for $ H^1_* \ii_X.$
\end{proposition}

\begin{proof} The proof of the first two statements rests on results from Section
\ref{properties-of-X}. The last statement follows from the
minimality of $ X $ in its biliaison class.
\end{proof}

\begin{proposition} $ h^0(X,
\oo_X(d)) = 4d+2b-1 $ for $ d \geq 2.$
\end{proposition}

\begin{proof} See Proposition \ref{h0OX}.
\end{proof}

As before, we describe the elements of $ H^0(X, \oo_X(d)) $ for
every $ d \geq 2.$ At first, we describe the elements of $ H^0(X,
\oo_X(d)) $ for $ d \geq b+1,$ because, in the considered range,
we have $ (R/I_X)_d = H^0(X, \oo_X(d)),$ and so it holds
\begin{proposition} Let $ d \geq b+1.$ Then, $ H^0(X, \oo_X(d)) $ is equal to
$$ \{ p_1 + y p_2 + (xz-y^2) p_3 + (xz-y^2) y p_4 + w
x^{d-b-1} p_5 + w x^{d-b-1} y p_6 \vert p_i \in K[x,z] \} $$ where
the degrees of the $ p_i'$ s are fixed in such a way that the
elements are homogeneous of degree $ d.$
\end{proposition}

To describe the elements of $ H^0(X, \oo_X(d)) $ for $ 2 \leq d
\leq b,$ this time we need two suitable global sections $ \xi_1,
\xi_2 $ of degree $ -b+2,$ and then we compute all the global
sections by using $ \xi_1, \xi_2 $ and the elements in $
(R/I_X)_d.$

\begin{definition} Let $ \xi_1, \xi_2 \in H^0(X, \oo_X(-b+2)) $ be the
global section of $ X $ defined as $$ \xi_1 = \frac{xz-y^2}{z^b} =
\frac{yw}{x^b} \mbox{ and } \xi_2 = \frac{y(xz-y^2)}{z^{b+1}} =
\frac{w}{x^{b-1}}.$$
\end{definition}

The global sections $ \xi_1 $ and $ \xi_2 $ are well defined because for no closed
point on $ C $ both $ x $ and $ z $ can be equal to $ 0,$ and
because the two definitions agree on the overlap (see the last two
generators of $ I_X$.)

\begin{proposition} \label{rel-xi1-xi2} $ \xi_1 $ and $ \xi_2 $  verify the following
equalities \begin{enumerate} \item $ w \xi_1 = (xz-y^2) \xi_1 = w
\xi_2 = (xz-y^2) \xi_2 = 0;$ \item $ x^b \xi_1 = y w, z^b \xi_1 =
xz-y^2;$ \item $ x^{b-1} \xi_2 = w, y z^b \xi_2 = x(xz-y^2),
z^{b+1} \xi_2 = y(xz-y^2);$ \item $ x \xi_1 = y \xi_2, y \xi_1 = z
\xi_2.$
\end{enumerate}
\end{proposition}

\begin{proof} The equalities easily follow from the definition of $ \xi_1, \xi_2 $ and
from the knowledge of the ideal $ I_X.$
\end{proof}

\begin{proposition} Let $ 2 \leq d
\leq b.$ Then $ H^0(X, \oo_X(d)) = V_d $ where \begin{equation*}
\begin{split} V_d = & \{p_1 + y p_2 + (xz-y^2) p_3 + (xz-y^2) y p_4 + w p_5
+ wy p_6 + \\ & x^{d-2} z^{d-1} \xi_2 (z q_1 + y q_2) \vert p_i
\in K[x,z], q_j \in K[x,z] \} \end{split}
\end{equation*} and the elements in $ V_d $ are homogeneous
of degree $ d.$
\end{proposition}

\begin{proof} We can apply the same argument as in Proposition
\ref{g-s-1}, with the only difference that $$ \{ f \in R \vert f
\xi_2 \in R/I_X \} = \langle w, xz-y^2, x^{b-1}, y z^b, z^{b+1}
\rangle $$ as can be easily computed via computer algebra
techniques.
\end{proof}

\begin{remark} \rm Thanks to the relations of previous Proposition \ref{rel-xi1-xi2}, we can write the elements in $ H^0(X, \oo_X(d)) $ without using $ \xi_1.$ Analogously, one can write them using $ \xi_1 $ but not $ \xi_2.$ The choice of using only one between $ \xi_1, \xi_2 $ allows to simplify the following computations.
\end{remark}

We want to compute the degree $ 0 $ global sections of the normal
sheaf of $ X $ as $$ H^0(X, \ns_X) = \ker( \Hom(F_1,
H^0_*(X,\oo_X)) \stackrel{\delta^{\vee}_1}{\lra} \Hom(F_2, H^0_*(X,
\oo_X)))_0 $$ where $ F_1 = R(-2) \oplus R(-3) \oplus R(-4) \oplus
R^2(-b-2) $ and $ F_2 = R(-4) \oplus R(-5) \oplus R^4(-b-3) $ and
$ \delta_1^{\vee} $ is the dual of $ \delta_1: F_2 \to F_1.$

To this end, let $ \varphi \in \Hom(F_1, H^0_*(X, \oo_X)) $ be a
degree $ 0 $ map. Then, if $ F_1 = \oplus_{i=1}^5 R e_i $ with $
\deg(e_1) = 2, \deg(e_2) = 3, \deg(e_3) = 4, \deg(e_4) = \deg(e_5)
= b+2,$ we have that $ \varphi(e_i) \in H^0(X, \oo_X(\deg(e_i))),$
and so, if we assume that $ b \geq 4,$ we can set
\begin{itemize} \item[$\bullet$] $ \varphi(e_1) = p_{11} + y p_{12} +
(xz-y^2) p_{13} + w p_{15} + w y p_{16} + \xi_2 (z^2 q_{11} + y z
q_{12}),$ with $ \deg(p_{11}) = 2, \deg(p_{12}) = \deg(p_{15}) =
1, \deg(p_{13}) = \deg(p_{16}) = 0, \deg(q_{11}) = \deg(q_{12}) =
b-2;$ \item[$\bullet$] $ \varphi(e_2) = p_{21} + y p_{22} +
(xz-y^2) p_{23} + (xz-y^2) y p_{24} + w p_{25} + w y p_{26} +
\xi_2 (x z^3 q_{21} + x y z^2 q_{22}),$ with $ \deg(p_{21}) = 3,
\deg(p_{22}) = \deg(p_{25}) = 2, \deg(p_{23}) = \deg(p_{26}) = 1,
\deg(p_{24}) = 0, \deg(q_{21}) = \deg(q_{22}) = b-3;$
\item[$\bullet$] $ \varphi(e_3) = p_{31} + y p_{32} + (xz-y^2)
p_{33} + (xz-y^2) y p_{34} + w p_{35} + w y p_{36} + \xi_2 (x^2
z^4 q_{31} + x^2 y z^3 q_{32}),$ with $ \deg(p_{31}) = 4,
\deg(p_{32}) = \deg(p_{35}) = 3, \deg(p_{33}) = \deg(p_{36}) = 2,
\deg(p_{34}) = 1, \deg(q_{31}) = \deg(q_{32}) = b-4;$
\item[$\bullet$] $ \varphi(e_4) = p_{41} + y p_{42} + (xz-y^2)
p_{43} + (xz-y^2) y p_{44} + w x p_{45} + w x y p_{46},$ with $
\deg(p_{41}) = b+2, \deg(p_{42}) = b+1, \deg(p_{43}) =
\deg(p_{45}) = b, \deg(p_{44}) = \deg(p_{46}) = b-1;$
\item[$\bullet$] $ \varphi(e_5) = p_{51} + y p_{52} + (xz-y^2)
p_{53} + (xz-y^2) y p_{54} + w x p_{55} + w x y p_{56},$ with $
\deg(p_{51}) = b+2, \deg(p_{52}) = b+1, \deg(p_{53}) =
\deg(p_{55}) = b, \deg(p_{54}) = \deg(p_{56}) = b-1.$
\end{itemize}

Of course, $ \varphi \in H^0(X, \ns_X) $ if, and only if, $
\delta_1^{\vee} (\varphi) = \varphi \circ \delta_1 = 0,$ and so we get
the following system \begin{equation} \label{system2} \left\{
\begin{array}{l} (xz-y^2) \varphi(e_1) - w \varphi(e_2) = 0 \\
(xz-y^2) \varphi(e_2) - w \varphi(e_3) = 0 \\ -y z^b \varphi(e_1)
+ x^b \varphi(e_2) - w \varphi(e_4) = 0 \\ z^b \varphi(e_2) - y
\varphi(e_4) + x \varphi(e_5) = 0 \\ -x^{b-1} \varphi(e_3) + z
\varphi(e_4) - y \varphi(e_5) = 0 \\ z^{b+1} \varphi(e_1) -
x^{b-1} y \varphi(e_2) + w \varphi(e_5) = 0. \end{array} \right.
\end{equation}

A technical result in the computation of $ H^0(X, \ns_X) $ is the
knowledge of the generators of two syzygy modules. In particular,
it holds
\begin{lemma} In $ R/I_X,$ the first syzygy module of $ (xz-y^2,
-w) $ is generated by $$ \left( \begin{array}{cccccc} xz-y^2 & w &
0 & 0 & x^b & x^{b-1} y \\ 0 & 0 & w & xz-y^2 & y z^b & z^{b+1}
\end{array} \right) $$ while the first syzygy module of $ (1, y) $
is minimally generated by $$ \left( \begin{array}{c} y \\ -1
\end{array} \right).$$
\end{lemma}

The proof is based on standard Gr\"obner bases arguments.

Thanks to the previous Lemma, we can solve system (\ref{system2}),
one equation at a time, and we get
\begin{theorem} \label{thm2} With the notation as above, $ h^0(X, \ns_X) = 9 + 4b = 13 - 2g.$
\end{theorem}

\begin{proof} To prove the statement, we have to compute the elements in $ H^0(X, \ns_X).$

\noindent \em Claim: \rm
$ \varphi \in H^0(X,
\ns_X) $ if, and only if, there exist $ \alpha, \beta \in K,$ and $ P_1,
\dots, P_8 \in K[x, z],$ of suitable degrees, such that
\begin{itemize} \item[$\bullet$] $ \varphi(e_1) = 2w (P_1 + y P_2);$
\item[$\bullet$] $ \varphi(e_2) = (xz-y^2) (P_1 + y P_2) + w (\alpha x^2 + z P_3 + y P_4);$
\item[$\bullet$] $ \varphi(e_3) = 2 (xz-y^2) (\alpha x^2 + z P_3 + y P_4);$
\item[$\bullet$] $ \varphi(e_4) = x^b
(\alpha x^2 + z P_3 + y P_4) - y z^b (P_1 + y P_2) + (xz-y^2) (x P_5 + y
P_6 + x^{b-1} P_3) + x w (x P_7 + y P_8 + \beta z^b);$ \item[$\bullet$]
$ \varphi(e_5) = x^{b-1} y (\alpha x^2 + z P_3 + y P_4) - z^{b+1} (P_1
+ y P_2) + (xz-y^2) (y P_5 + z P_6 - x^{b-1} P_4) + x w (y
P_7 + z P_8 - \alpha z^b).$ \end{itemize}

If the claim holds, then we can compute $ h^0(X, \ns_X) $ with an easy parameter count. Hence, we prove the claim. Its proof is a quite long computation where we use the same ideas
as in the proof of Theorem \ref{thm1}.
\end{proof}

\begin{remark} \rm Now, we consider the cases not covered by
Theorem \ref{thm2}, namely $ 1 \leq b \leq 3.$ We consider the double conic
defined by the ideal $$ I_X = \langle w^2, w(xz-y^2), (xz-y^2)^2, x^b
(xz-y^2) - y z^b w, x^{b-1} y (xz-y^2) - z^{b+1} w \rangle $$
with $ b = 1, 2, 3,$ and we compute $ h^0(X, \ns_X) $ by using Macaulay (see \cite{mac}).

If $ b = 3,$ we get $ h^0(X, \ns_X) = 21 = 13 - 2 g,$ because $ g = 2 - 2b = -4.$

If $ b = 2,$ we get $ h^0(X, \ns_X) = 17 = 13 - 2 g,$ because $ g = -2.$

If $ b = 1,$ then we get $ H^0(X, \ns_X) = 16 \not= 13 - 2g,$ because, in this case, $ g = 0.$
\end{remark}

\subsection{Case $ X \subseteq \pp^n, n \geq 4$}

Now, we suppose that $ X $ is a suitable double conic in $ \pp^n, n \geq 4.$
\begin{theorem} \label{thm3} Let $ C \subset \pp^n $ be the conic defined by the ideal
$$ I_C = \langle x_0 x_2 - x_1^2, x_3, \dots, x_n \rangle,$$ and let $ j: \pp^1 \to C $
be the isomorphism defined as $ j(t:u) = (t^2: t u: u^2: 0: \dots
: 0).$ Let $ \mu: \oo_{\pp^1}(-4) \oplus \oo_{\pp^1}^{n-2}(-2) \to
\oo_{\pp^1}(-4 + a), a \geq 5,$ be the map defined as $ \mu =
(u^a, t^{a-2}, 0, \dots, 0),$ and let $ X $ be the double
structure on $ C $ associated to $ \mu $ and $ j.$ Then, $ h^0(X,
\ns_X) \leq (n-1) (5-g) + 3.$
\end{theorem}

\begin{proof} $ \oo_{\pp^1}^{n-2}(-2) $ is contained in the kernel of $ \mu.$ Hence, $ x_3, \dots, x_n \in I_X $ i.e. $ X $ is degenerate and it is contained in the linear space $ L $ of dimension $ 3.$ Then, we
can consider both the normal sheaf $ \ns_{X,L} $ of $ X $ in $ L,$ and the normal sheaf
$ \ns_X  $ of $ X $ in $ \pp^n.$ They are related via the exact sequence $$ 0 \to \ns_{X,L}
\lra \ns_X  \lra \oo_X(1)^{n-3}.$$ In particular, we have the inequality $$ h^0(X, \ns_X)
\leq h^0(X, \ns_{X,L}) + (n-3) h^0(X, \oo_X(1)).$$ By hypothesis, $ a \geq 5 $ and so the
arithmetic genus $ g $ of $ X $ satisfies $ g \leq -2.$ By Theorems \ref{thm1}, \ref{thm2}
and the Remarks after them, $ h^0(X, \ns_{X,L}) = 13 - 2g,$ while $ h^0(X, \oo_X(1)) =
\dim_K (R/I_X)_1 + h^1 \ii_X(1) = 4 + a-2 = 5-g $ as proved in Propositions \ref{rao-f}. Then,
$$ h^0(X, \ns_X) \leq  13 - 2g + (n-3) (5-g) = (n-1) (5-g) + 3.$$
\end{proof}

\subsection{Remarks on the Hilbert schemes $ \Hilb_{4t+1-g}(\pp^n) $}

In this last subsection, we use the previous results to get
information on the irreducible components
containing the double structures of genus $ g $ on conics.

\begin{theorem} If $ g \leq -2,$ then $ \overline{H(4,g,n)} $ is a generically
smooth irreducible component of $ \mathcal{H}ilb_{4t+1-g}(\pp^n) $ of dimension $ (n-1)(5-g) + 3.$
\end{theorem}

\begin{proof} As proved in Corollary \ref{dim-comp}, every irreducible component containing $ \overline{H(4,g,n)} $ has
dimension greater than or equal to $ (n-1)(5-g) + 3,$ but the tangent space to the Hilbert scheme $ \mathcal{H}ilb_{4t+1-g}(\pp^n) $ at the double conic described in Theorem
\ref{thm3} has dimension lesser than or equal to $ (n-1)(5-g) + 3.$ Hence, the point corresponding
to the double conic considered in Theorem \ref{thm3} is smooth, $ \overline{H(4,g,n)} $ is irreducible of dimension $ \dim \overline{H(4,g,n)} = (n-1)(5-g) + 3 $ and the point corresponding
to the double conic considered in Theorem \ref{thm3} is smooth,
i.e. $ \overline{H(4,g,n)} $ is generically smooth.
\end{proof}

Now, we add some remarks to $ \overline{H(4,g,3)} $ for $ g \geq
-1.$

\begin{proposition} If $ g = -1,$ the general element of $ \overline{H(4,-1,3)} $ is
the union of two smooth conics without common points, and a double
structure on a smooth conic is a smooth point of $
\overline{H(4,-1,3)}.$
\end{proposition}

\begin{proof} By a simple parameter count, the family of two disjoint
conics has dimension $ 16,$ which is equal to the dimension of the tangent
space to $ \overline{H(4,-1,3)} $ at the double conic considered
in Theorem \ref{thm1} and in the subsequent Remark. The claim
follows if we exhibit a family whose general element is a disjoint
union of two conics, and whose special fiber is the considered
double structure. The ideal $$ \langle w, xz-y^2 \rangle \cap
\langle w+tx, tz^2 + xz - y^2 \rangle \subseteq K[x,y,z,w,t] $$
gives a flat family over $ \aa^1 $ with the required properties.
\end{proof}

\begin{proposition} If $ g = 0,$ the general element of $ \overline{H(4,0,3)} $ is
a rational quartic curve, and a double structure on a smooth conic
is a smooth point of $ \overline{H(4,0,3)}.$
\end{proposition}

\begin{proof} By a simple parameter count, the family of the rational
quartic curves has dimension $ 16,$ which is equal to the dimension of the
tangent space to $ \overline{H(4,0,3)} $ at the double conic
considered in Theorem \ref{thm2} and in the subsequent Remark. The
claim follows if we exhibit a family whose general element is a
rational quartic curve, and whose special fiber is the considered
double structure. The ideal $$ \langle w^2+t(xy-zw),
y(xz-y^2)-z^2w \rangle : \langle xy-zw, y^2, yw, w^2 \rangle
\subseteq K[x,y,z,w,t] $$ gives a flat family over $ \aa^1 $ with
the required properties. In fact, a general quartic curve is
linked to two skew lines or to a double line of genus $ -1 $ via a
complete intersection of type $ (2,3).$ The general element of the
family is the residual intersection of a double line of genus $ -1
$ on a smooth quadric surface, while the special element is the
residual intersection to the same double line on a double plane.
\end{proof}

\begin{proposition} If $ g = 1,$ the general element of $ \overline{H(4,1,3)} $
is the complete intersection of two quadric surfaces, and a double
structure on a smooth conic is a smooth point of $
\overline{H(4,1,3)}.$
\end{proposition}

\begin{proof} By a simple parameter count, the family of the complete
intersections of two quadric surfaces  has dimension $ 16,$ which is equal
to the dimension of the tangent space to $ \overline{H(4,1,3)} $
at the double conic considered in Theorem \ref{thm1} and in the
subsequent Remark. The claim follows because it is easy to check
that the ideal of the considered double conic is the complete
intersection of $ w^2 $ and $ xz-y^2-zw.$
\end{proof}

\begin{proposition} If $ g = 3,$ the general element of $ \overline{H(4,3,3)} $ is
a plane quartic curve, and a double structure on a smooth conic is
a smooth point of $ \overline{H(4,3,3)}.$
\end{proposition}

\begin{proof} By a simple parameter count, the family of the plane
quartic curves has dimension $ 17,$ which is equal to the dimension of the
tangent space to $ \overline{H(4,3,3)} $ at the double conic
considered in Theorem \ref{thm1} and in the subsequent Remark. The
claim follows because the double conic we considered is a plane
quartic curve.
\end{proof}



\begin{thebibliography}{11}

\bibitem{bbm} E. Ballico, G. Bolondi, J.C. Migliore, {\em The Lazarsfeld-Rao problem for liaison classes of two-codimensional subschemes of $ \pp^n $}, Amer. J. Math. \bf 113 \rm (1991),  no. 1, 117--128.

\bibitem{mac} D. Bayer, M. Stillmann, {\em Macaulay: a system for
computation in algebraic geometry and commutative algebra,} Source
and object code available for Unix and Macintosh computers.
Contact the authors, or download from zariski.harvard.edu via
anonymous ftp. (login: anonymous, password: any, cd Macaulay),
Cambridge U.P., 1982-1990.

\bibitem{bd} D. Bayer, D. Eisenbud, {\em Ribbons and their canonical
embeddings}, Trans. Amer. Math. Soc. \bf 347 \rm (1995), no. 3, 719--756.

\bibitem{bafo} C. Banica, O. Forster, {\em Multiplicity structures
on space curves}, Contemporary Math. \bf 58 \rm (1986), 47--64.

\bibitem{bama} C. Banica, N. Manolache, {\em Rank $ 2 $ stable vector bundles on $ \pp^3(\mathbb C) $ with Chern classes $ c_1 = -1, c_2 = 4 $}, Math. Z. \bf 190 \rm (1985), 315-339.

\bibitem{bg} M. Boraty\'nski, S. Greco, {\em When does an ideal
arise from the Ferrand construction?}, Boll. Un. Mat. Ital. B (7) 1 (1987), no. 1, 247--258.

\bibitem{greco-bor-2} M. Boraty\'nski, S. Greco,{\em Hilbert functions and Betti numbers in a flat family}, Ann. Mat. Pura Appl. (4)  142  (1985), 277--292 (1986).

\bibitem{cm} M. Casanellas, R.M. Mir\'o-Roig, {\em On the Lazarsfeld-Rao
property for Gorenstein liaison classes},  J.  Pure and Appl. Alg.
\bf 179 \rm (2003), 7-12.

\bibitem{eisenbud-book} D. Eisenbud, Commutative algebra, \bf GTM 150, \rm Springer-Verlag, 2004.

\bibitem{eg} D. Eisenbud, M. Green, {\em Clifford indices of
ribbons}, Trans. Amer. Math. Soc. \bf 347 \rm (1995), no. 3, 757--765.

\bibitem{fer} D. Ferrand, {\em Courbes gauches et fibres de rang
2}, C.R. Acad. Sci. Paris Ser. A \bf 281 \rm (1977), 345--347.

\bibitem{fong} L.Y. Fong, {\em Rational ribbons and deformation
of hyperelliptic curves}, J. Algebraic Geom. \bf 2 \rm (1993), no. 2, 295--307.

\bibitem{g} F. Gaeta, {\em Nuove ricerche sulle curve sghembe algebriche
di risiduale finito e sui gruppi di punti del piano}, (French) An. Mat. Pura Appl. (4) 31, (1950), 1--64.

\bibitem{gr} A. Grothendieck ,{\em Techniques de construction et th\'eor\`{e}mes d'
existence en g\'{e}om\'{e}trie alg\'{e}brique. IV. Le sch\'{e}mas
de Hilbert}, S\'{e}minaire Bourbaki, Exp. No. 221, vol. 8,
Soc.Math.France, Paris, 1995, pp.249-276.

\bibitem{hart1} R. Hartshorne, {\em Connectedness of the Hilbert scheme}, Publ.
Math. de I.H.E.S. \bf 29 \rm (1966), pp. 261-304.

\bibitem{hart2} R. Hartshorne, {\em Some examples of Gorenstein
liaison in codimension three}, Collect. Math. \bf 53 \rm (2002),
21-48.

\bibitem{hart} R. Hartshorne, Algebraic Geometry, {\bf GTM 52},
Springer-Verlag, 1977.

\bibitem{hs} R. Hartshorne, E. Schlesinger, {\em Curves in the double
plane}, Special issue in honor of Robin Hartshorne, Comm. Algebra \bf 28 \rm (2000), no. 12, 5655--5676.

\bibitem{kempf-book} G.R. Kempf, {\em Algebraic varieties}, London Mathematical Society Lecture Note Series, \bf 172\rm, Cambridge University Press, Cambridge, 1993.

\bibitem{man1} N. Manolache, {\em Multiple structures on smooth support}, Math. Nachr. \bf 167 \rm (1994), 157--202.

\bibitem{man2} N. Manolache, {\em Double rational normal curves with
linear syzygies}, Manuscripta Math. \bf 104 \rm (2001), no. 4, 503--517.

\bibitem{mdp} M. Martin-Deschamps, D. Perrin, {\em Sur la
classification des courbes gauches}, Ast\'erisque 184-185, Soc.
Math. de France (1990).

\bibitem{mig-book} J.C. Migliore, Introduction to liaison theory
and deficiency modules, Progress in Mathematics \bf 165 \rm
Birkh\"auser, 1998.

\bibitem{nns2} U. Nagel, R. Notari, M.L. Spreafico, {\em On the even
Gorenstein liaison classes of ropes on a line}, Dedicated to Silvio
Greco on the occasion of his 60th birthday (Catania, 2001). Matematiche
(Catania) 55 (2000), no. 2, 483--498 (2002).

\bibitem{nns1} U. Nagel, R. Notari, M.L. Spreafico, {\em Curves of degree
two and ropes on a line: their ideals and even liaison classes}, J.
Algebra \bf 265 \rm (2003), no. 2, 772--793.

\bibitem{nns3} U. Nagel, R. Notari, M.L. Spreafico, {\em The Hilbert scheme
of degree two curves and certain ropes}, Internat. J. Math. \bf 17 \rm (2006), no. 7, 835--867.

\bibitem{ns} S. Nollet, E. Schlesinger, {\em Hilbert schemes of degree
four curves}, Compositio Math. \bf 139 \rm (2003), no. 2, 169--196.

\bibitem{sernesi} E. Sernesi, Deformations of algebraic schemes, Grundlehren der Mathematischen Wissenschaften \bf 334, \rm Springer-Verlag, Berlin, 2006.

\bibitem{wahl} J. Wahl, {\em  On cohomology of the square of an ideal sheaf},
J. Algebraic Geom. \bf 6 \rm (1997), no. 3, 481--511.

\end{thebibliography}
\end{document}